\documentclass[onefignum,onetabnum]{siamart171218}
\usepackage{geometry}
\usepackage[utf8]{inputenc}
\usepackage{xcolor}
\usepackage[english]{babel}
\usepackage{amsfonts,amssymb,amstext}
\usepackage{mathtools}
\usepackage[shortlabels]{enumitem}
\usepackage[tight,nice]{units}
\usepackage{url}
\usepackage{graphicx}
\usepackage{bbm}
\usepackage[parse-numbers=false,per-mode=reciprocal]{siunitx}
\usepackage{textcomp}
\usepackage{algorithm}
\usepackage{algpseudocode}
\usepackage{multirow}
\usepackage{booktabs}

\newtheorem{thm}{Theorem}
\numberwithin{thm}{section}

\newtheorem{lem}[thm]{Lemma}

\newtheorem{rem}[thm]{Remark}
\newtheorem{ass}[thm]{Assumption}

\newcommand{\bu}{\mbox{\boldmath$u$}}

\newcommand{\bU}{\mbox{\boldmath$U$}}
\newcommand{\bj}{\mbox{\boldmath$j$}}
\newcommand{\bz}{\mbox{\boldmath$z$}}

\newcommand{\IKR}[1]{\textcolor{black}{#1}} %
\newcommand{\IKRr}[1]{\textcolor{black}{#1}} %

\newcommand{\IKRrr}[1]{\textcolor{black}{#1}} 
\newcommand{\todo}[1]{\textcolor{black}{#1}}

\headers{Parallel-in-Time Multi-Harmonic Solution of Time-Periodic Problems}{Kulchytska-Ruchka and Schöps}

\title{Efficient Parallel-in-Time Solution of Time-Periodic Problems Using a Multi-Harmonic Coarse Grid Correction}

\author{%
	Iryna Kulchytska-Ruchka\thanks{
		Computational Electromagnetics Group, Technical University of Darmstadt, Schlossgartenstra{\ss}e 8, D-64289 Darmstadt, Germany.}
	\and
	Sebastian Sch\"ops\footnotemark[1]}

\begin{document}

	\maketitle
	
	\begin{abstract}
	This paper presents a highly-parallelizable parallel-in-time algorithm 
	for efficient solution of nonlinear time-periodic problems. It is 
	based on the time-periodic extension of the Parareal method, known to 
	accelerate sequential computations via parallelization on the fine grid. 
	The proposed approach reduces the complexity of the periodic Parareal solution 
	by introducing a simplified Newton algorithm, which allows an 
	additional parallelization on the coarse grid. In particular, at each Newton 
	iteration a multi-harmonic correction is performed, which converts the 
	block-cyclic periodic system in the time domain into a block-diagonal system 
	in 	the frequency domain, thereby 
	\IKRrr{solving} for each frequency 
	component in parallel. \IKRrr{The convergence analysis of the method is discussed for a one-dimensional model problem.} 
	The introduced algorithm and several existing solution 
	approaches are compared via their application to the eddy current problem for 
	both linear and nonlinear models of a coaxial cable. Performance of the 
	considered methods is also illustrated for a three-dimensional transformer model.
	\end{abstract}
	
	

	\section{Introduction} \label{section:introduction}	
	{The computation of} the time-periodic solution of a \IKR{stiff} evolution problem becomes 
	particularly \IKR{time consuming}, when one is interested in the steady-state 
	behavior 
	e.g., \IKR{of} an electrical machine. 
	\IKR{The} numerical treatment of time-periodic problems is often {challenging} 
	since it requires \IKR{the} simultaneous solution of the underlying system on 
	the whole period. {Indeed, a classical methodology involves a 
	discretization in space, e.g., using the finite element method (FEM), followed by a finite difference discretization 
	in time, which couples all the discrete values on one period through the periodicity condition. 
	This space-time discretization approach of time-periodic problems is known in engineering as 
	the time-periodic (TP) FEM \cite{Hara_1985aa}. 
	It} leads to a large system of algebraic equations  
	whose solution would usually be very computationally expensive 
	or even prohibitive, especially due to 
	its troublesome block-cyclic structure. 
	These obstacles on the way to solve time-periodic problems cause 
	an urgent need to develop efficient \IKR{parallel} algorithms, 
	able to simplify and accelerate the time-domain computations.
	
	Based on the idea of the multiple shooting method \cite{Stoer_2005aa}, the 
	Parareal algorithm was originally introduced to speed up sequential solution 
	of initial value problems (IVPs) using parallelization in the time domain 
	\cite{Lions_2001aa}. It involves computations on two grids (fine and 
	coarse). While the coarse solver gives rough information about the 
	solution and is applied sequentially, the acceleration is obtained via 
	parallel calculation of the accurate fine solution. Convergence properties 
	of the method were investigated in \cite{Gander_2007aa}, \cite{Gander_2008aa}. 
	Application of Parareal to the simulation of an electrical machine in \cite{Schops_2018aa} 
	illustrated its significant acceleration capabilities. The original Parareal 
	approach for IVPs was subsequently extended to the class of time-periodic 
	problems and analyzed in \cite{Gander_2013ab}. 

	This work introduces a novel Parareal-based iterative algorithm able to solve nonlinear time-periodic 
	problems efficiently. In contrast to the fixed point iteration introduced in \cite{Gander_2013ab} 
	a special linearization via the simplified Newton method is proposed here. A specific choice of 
	the initial approximation for Newton's iteration gives a block-cyclic system matrix, which can be 
	diagonalized using the frequency domain approach from \cite{Biro_2006aa}. Following the terminology of  
	\cite{Bachinger_2005aa} we call it the multi-harmonic (MH) solver in this paper. This technique was recently   
	used for parallel-in-time solution of linear time-periodic problems in \cite{Kulchytska-Ruchka_2019ac}. 
	It was applied on the coarse level of the time-periodic Parareal system, 
	\IKRrr{calculating} each frequency component separately and in parallel. Diagonalization using the fast Fourier transform (FFT)
	\cite{Trefethen_1996aa} makes the MH approach {particularly} attractive due to its low complexity. \IKRrr{The Fourier-based diagonalization was exploited to find a time-periodic solution of a fractional diffusion equation in \cite{Wu_2018aa}, where each separate system of equations was solved using the multigrid method. Besides,} a similar idea of diagonalization and parallelization was proposed for solution of IVPs 
	in \cite{Gander_2017aa}, \cite{Gander_2019ab}, 
	where non-equidistant partition of the considered time interval was required. \IKRrr{Finally, 
	the MH approach was recently applied to the parallel-in-time solution of IVPs in \cite{Gander_2020aa}. 
	There an additional restructuring to a block-cyclic form was performed.}
		
	
	Our paper is organized as follows. At first, a couple of existing parallel-in-time approaches 
	for time-periodic problems \cite{Gander_2013ab} are described in Section~\ref{section:PP-PC}. 
	Section~\ref{section:Parareal_HB} recalls the main ideas of the frequency domain solution \IKRrr{of linear time-periodic problems. In Section~\ref{section:MH_TP} the MH approach presented} in \cite{Biro_2006aa} as well as its parallel-in-time extension \cite{Kulchytska-Ruchka_2019ac} \IKRrr{are described}. A special \IKR{Newton-type solver} is then introduced in 
	Section~\ref{section:Parareal_new} for solution of nonlinear time-periodic problems with 
	the MH approach. \IKR{Its convergence properties are analyzed using the theory from \cite{Deuflhard_2004aa} 
	and illustrated for a scalar model problem.} 
	Application of the proposed method to numerical treatment of the eddy 
	current problem \IKRr{for a coaxial cable model}, as well as its performance compared to the conventional sequential time stepping,  
	to the time-parallel algorithms from \cite{Gander_2013ab}, and to the MH solution of  
	\cite{Biro_2006aa} are presented in Section~\ref{sec:numerics}. \IKRr{Section~\ref{sec:numerics_transfo} describes
	performance of the considered approaches applied to a three-dimensional model of a transformer.} 
	The paper is finally concluded with Section~\ref{section:conclusions}.	
	\section{\IKR{Overview of the }periodic Parareal-based parallel-in-time algorithms} \label{section:PP-PC}	
	Consider a time-periodic problem for a system of ordinary differential equations (ODEs) 
	on time interval $(0,T)$ 
\begin{equation} \label{eq:ode}
\begin{aligned}
\mathbf{M}\bu^{\prime}(t) + \mathbf{K}\bigl(\bu(t)\bigr)\bu(t)&=\bj(t), \quad t \in (0,T),\\
  \qquad\bu(0) &= \bu(T),
  \end{aligned}
\end{equation}
with unknown $\bu:[0,T]\to\mathbb{R}^d,$ matrices $\mathbf{M}$ and $\mathbf{K}(\bu(t))$ such that \IKRrr{the matrix pencil $\bigl(\mathbf{M},\mathbf{K}(\mathbf{v})\bigr)=\mathbf{K}(\mathbf{v})+\lambda\mathbf{M},$ with $\lambda\in\mathbb{R}$ is regular for each $\mathbf{v}\in\mathbb{R}^d,$}  
and a $T$-periodic right-hand side (RHS) $\bj:[0,T]\to\mathbb{R}^d$. 
System \eqref{eq:ode} could stem, e.g., from a spatial discretization 
of a parabolic partial differential equation (PDE) discretized by FEM \cite{Brenner_2008aa}.

As in the classical Parareal method \cite{Lions_2001aa}, we initially split the time interval $[0,T]$ 
into $N$ subintervals $[T_{n-1},T_n],$ $n=1,\dots,N$ using partition $0=T_0<T_1<\dots<T_N=T.$
An IVP for unknown $\bu_n:[T_{n-1}, T_{n}]\to\mathbb{R}^d$ given by 
\begin{equation}
\label{eq:ode_nf}
\begin{aligned}
\mathbf{M}\bu_n^{\prime}(t) + \mathbf{K}\bigl(\bu_n(t)\bigr)\bu_n(t)&=\bj(t),\quad t\in(T_{n-1}, T_{n}],\\
    \bu_n(T_{n-1}) &= \bU_{n-1}
\end{aligned}
\end{equation}
is considered on \IKR{the} $n$th subinterval, $n=1,\dots,N$. Within the Parareal setting, 
two propagators are applied to the {$n$th} IVP: the fine $\mathcal{F}\bigl(t, T_{n-1},{\mathbf{U}_{n-1}}\bigr)$ 
and the coarse $\mathcal{G}\bigl(t, T_{n-1},{\mathbf{U}_{n-1}}\bigr),$ \cite{Gander_2007aa}. Both 
operators calculate a solution to \eqref{eq:ode_nf} at $t\in(T_{n-1}, T_{n}],$ starting from \IKR{the} initial value 
$\bU_{n-1}$ at $T_{n-1},$ $n=1,\dots,N.$ However, accuracies and therefore computational costs 
of the two solvers differ. In particular, $\mathcal{F}$ solves the IVP with a very high precision, e.g., 
via time stepping over a very fine grid, while $\mathcal{G}$ provides only a rough approximation of the 
solution, e.g., by using a very coarse discretization or a lower-order method, compared to the fine propagator. 

	Based on the idea of Parareal, \IKR{Gander et al.} \cite{Gander_2013ab} introduce two parallel-in-time algorithms for efficient treatment of time-periodic problems. The first one, called periodic Parareal algorithm with initial-value coarse problem (PP-IC), calculates the periodic solution using \IKR{the} iteration
\begin{align}
\mathbf{U}_0^{(k+1)}&
=\mathbf{U}_N^{(k)},\label{eq:PP-IC1}\\
{{\mathbf{U}_n^{(k+1)}}}&
=
{\mathcal{F}}\bigl(T_n, T_{n-1},{\mathbf{U}^{(k)}_{n-1}}\bigr)+
{{\mathcal{G}}}\bigl(T_n, T_{n-1},{\mathbf{U}^{(k+1)}_{n-1}}\bigr)-{{\mathcal{G}}}\bigl(T_n, T_{n-1},{\mathbf{U}^{(k)}_{n-1}}\bigr)\label{eq:PP-IC2}
\end{align}
for $k=0,1,\dots$ and $n=1,\dots,N.$ While in the classical Parareal method the initial value (at $T_0$) is fixed, 
\eqref{eq:PP-IC1}-\eqref{eq:PP-IC2} updates the solution at $T_0$ 
with the solution at $T_N,$ obtained from the previous iteration. In this way PP-IC weakens the periodic 
coupling among the discrete values at the synchronization points $T_n,$ $n=0,\dots,N-1$ and, like Parareal, \IKR{it}
involves \IKR{the} solution of \IKRrr{a sequence of} IVPs \IKRrr{(and not of a single time-periodic problem)} on the coarse level. The method has been recently applied to simulation of an 
induction motor for an electric vehicle drive by the authors in \cite{Bast_2019aa} and already delivered 
a significant speedup. However, since PP-IC imposes a relaxed periodicity constraint one \IKR{may} 
expect a rather slow convergence to the periodic solution, especially when the underlying dynamical 
system possesses a very long settling time. \IKR{This is observed in the original paper \cite{Gander_2013ab} 
as well as within our numerical experiments in Section~\ref{sec:numerics}.}
	
The second approach presented in \cite{Gander_2013ab} is PP-PC: periodic Parareal algorithm with periodic coarse problem. 
In contrast to PP-IC, it maintains the prescribed periodic coupling between the first and the last values of the period 
$[T_0,T_N]$. The PP-PC iterations are written for $k=0,1,\dots$ as follows: 
\begin{align}
\mathbf{U}_0^{(k+1)}&
=\mathcal{F}\bigl(T_N, T_{N-1},\mathbf{U}^{(k)}_{N-1}\bigr)+\mathcal{G}\bigl(T_N, T_{N-1},\mathbf{U}^{(k+1)}_{N-1}\bigr) - \mathcal{G}\bigl(T_N, T_{N-1},\mathbf{U}^{(k)}_{N-1}\bigr),\label{eq:PP-PC1}\\
\mathbf{U}_n^{(k+1)}&
=
\mathcal{F}\bigl(T_n, T_{n-1},\mathbf{U}^{(k)}_{n-1}\bigr)+\mathcal{G}\bigl(T_n, T_{n-1},\mathbf{U}^{(k+1)}_{n-1}\bigr) - \mathcal{G}\bigl(T_n, T_{n-1},\mathbf{U}^{(k)}_{n-1}\bigr)\label{eq:PP-PC2}
\end{align}
with $n=1,\dots,N-1.$ It can be seen that the fine solutions only involve values from 
iteration $k$ and can therefore be computed in parallel, in the same way as within 
Parareal and PP-IC. On the other hand, the coarse propagation can \IKR{neither} be  
performed sequentially \IKR{nor in parallel}, since the periodicity condition yields interconnection of 
the solution at \IKR{all} synchronization points. Hence, PP-PC requires \IKR{the} solution of a 
periodic coarse problem, which implies a joint computation of the coarse grid 
values over the whole period. 

Denoting the difference between the fine and the coarse solutions of IVP \eqref{eq:ode_nf} at $T_n,$ 
calculated starting from \IKR{the} initial value $\mathbf{U}^{(k)}_{n-1}$ at $T_{n-1}$, by
\begin{equation}
\label{eq:bnk}
\mathbf{b}_n^{(k)}=\mathcal{F}\bigl(T_n, T_{n-1},\mathbf{U}^{(k)}_{n-1}\bigr)-\mathcal{G}\bigl(T_n, T_{n-1},\mathbf{U}^{(k)}_{n-1}\bigr),\quad n=1,\dots,N 
\end{equation}
we write PP-PC \IKR{the} iteration \eqref{eq:PP-PC1}-\eqref{eq:PP-PC2} in a matrix-vector (operator) form 
\begin{equation}
\label{eq:PP-PC_sys}
\setlength{\arraycolsep}{2.5pt}
\begin{bmatrix}
\mathbf{I} & \mathbf{0}   & \dots & -\mathcal{G}\left(T_N,T_{N-1},\cdot\right) \\
-\mathcal{G}\left(T_1,T_0,\cdot\right)  & \mathbf{I} &  & \mathbf{0} \\
\vdots &  \ddots   &  \ddots &\vdots \\
\mathbf{0}&    \dots &  -\mathcal{G}\left(T_{N-1},T_{N-2},\cdot\right) & \mathbf{I}
\end{bmatrix}
	\begin{bmatrix}
	\mathbf{U}_0^{(k+1)}\\
	\mathbf{U}_1^{(k+1)}\\
	\vdots\\
	\mathbf{U}_{N-1}^{(k+1)}
	\end{bmatrix}=\begin{bmatrix}
		\mathbf{b}_N^{(k)}\\
		\mathbf{b}_1^{(k)}\\
		\vdots\\
		\mathbf{b}_{N-1}^{(k)}
\end{bmatrix},
\end{equation}
where $\mathbf{I}\in\mathbb{R}^{d\times d}$ denotes the identity matrix. Note that the system matrix of \eqref{eq:PP-PC_sys} might have large dimensions, 
since each of the unknown $\mathbf{U}_n^{(k+1)},$ $n=0,\dots,N-1$ possibly consists of many degrees of freedom, e.g., coming from a spatial discretization. 
The matrix also possesses a block-cyclic-like structure, which is typical for time-periodic systems but inconvenient for many linear solvers. 

System \eqref{eq:PP-PC_sys} is written in an implicit operator form and cannot be solved directly in the nonlinear case.
A Jacobi-type fixed point iteration was applied in \cite{Gander_2013ab} to iteratively solve the PP-PC system. It is 
written for $s=0,1,\dots$ as
\begin{equation}
\label{eq:PP-PC_sys_fxd_pt}
	\begin{bmatrix}
	\mathbf{U}_0^{(k+1,s+1)}\\
	\mathbf{U}_1^{(k+1,s+1)}\\
	\vdots\\
	\mathbf{U}_{N-1}^{(k+1,s+1)}
	\end{bmatrix}=
	\setlength{\arraycolsep}{2.5pt}
	\begin{bmatrix}
\mathbf{0} & \mathbf{0}   & \dots & \mathcal{G}\left(T_N,T_{N-1},\cdot\right) \\
\mathcal{G}\left(T_1,T_0,\cdot\right)  & \mathbf{0} &  & \mathbf{0} \\
\vdots &  \ddots   &  \ddots &\vdots \\
\mathbf{0}&    \dots &  \mathcal{G}\left(T_{N-1},T_{N-2},\cdot\right) & \mathbf{0}
\end{bmatrix}
	\begin{bmatrix}
	\mathbf{U}_0^{(k+1,s)}\\
	\mathbf{U}_1^{(k+1,s)}\\
	\vdots\\
	\mathbf{U}_{N-1}^{(k+1,s)}
	\end{bmatrix}+
	\begin{bmatrix}
\mathbf{b}_N^{(k)}\\
\mathbf{b}_1^{(k)}\\
\vdots\\
\mathbf{b}_{N-1}^{(k)}
\end{bmatrix}
\end{equation}
{at each} PP-PC iteration {$k+1$.} {The} fixed point iteration \eqref{eq:PP-PC_sys_fxd_pt} decouples \IKR{the} calculation of the values $\mathbf{U}_n^{(k+1,s+1)},$ $n=0,\dots,N-1$ 
completely, thereby allowing parallel solution even on the coarse level. \IKRrr{Another advantage of this iterative method is that the  is fully non-intrusive, i.e., the propagators $\mathcal{F}$ and $\mathcal{G}$ can be considered as black-box solvers.} On the other hand, such a lineariazation relaxes the periodicity constraint similarly to 
PP-IC and \IKR{thus a loss of efficiency may be expected. This becomes obvious in the linear case where \eqref{eq:PP-PC_sys_fxd_pt} is equivalent to the block Jacobi method for linear systems.}

Our aim is to construct an iterative method, which benefits from the present block-cyclic structure and provides fast convergence. 

\section{\IKR{MH frequency domain solution of time-periodic problems}}\label{section:Parareal_HB}
\IKR{The} frequency domain representation \IKR{of \eqref{eq:ode}} is often exploited in (electrical) engineering \IKRrr{since it transforms} 
a differential equation into an algebraic one \cite{Nakhla_1976aa}. 
In contrast to the possibly lengthy time stepping in the time domain, MH frequency domain computations allow to obtain the steady-state solution 
directly via calculating its harmonic components. However, it may require many basis functions if the solution exhibits local features, e.g., due to pulsed excitations \cite{Bose_2006aa, Mohan_2003aa}. 
This frequency domain solution approach is also called the harmonic balance method 
\cite{Gyselinck_2002aa}, \cite{Auserhofer_2007aa} and  
can be interpreted as the Fourier collocation \cite{Deuflhard_2004aa} or the Fourier spectral method \cite{Trefethen_1996aa}. 

\IKRrr{Let $u$ be} a periodic function from the Lebesgue space 
$L^2([0,T],\mathbb{C}),$ \IKRrr{i.e., $u:[0,T]\to\mathbb{C}$ is a measurable function with a square integrable norm $\|u(t)\|_{\mathbb{C}}$ on $[0,T]$ (see \cite[Chapter 23.2]{Zeidler_1990aa}) and $u(0)=u(T)$.} 
The main idea of the MH solution \cite{Biro_2006aa} \IKRrr{is to represent $u$}
\IKRrr{with} the infinite sum 
$$
u(t)=\sum\limits_{j=-\infty}^{\infty}\hat{u}_j\psi_j(t),\quad{t\in[0,T]}
$$ 
of the spectral (orthonormal) basis functions $\psi_j(t)=e^{\imath {{2\pi}jt/{T}}}$ on time interval $[0,T],$ 
with (Fourier) coefficients
\begin{equation}
\label{eq:FT}
\hat{u}_j:=\frac{1}{T}\int\limits_0^T u(t)\psi_j^*(t)dt=\frac{1}{T}\int\limits_0^T u(t)e^{-\imath {{2\pi}jt/{T}}}dt,\quad j\in\mathbb{Z}.
\end{equation} 
Searching for the $T$-periodic solution $u$ 
of an ODE
\begin{equation}\label{eq:HB_ode}
F(u^{\prime},u,t)=0,\quad t\in(0,T), 
\end{equation}
we write its variational formulation
\begin{equation}
\label{eq:var_form}
0=(F(u^{\prime},u,t),\psi)_{L^2([0,T],\mathbb{C})}:=\frac{1}{T}\int\limits_0^T F(u^{\prime},u,t)\psi^*(t)dt\quad\forall\psi\in L^2([0,T],\mathbb{C}),
\end{equation} 
{where it is assumed that $u^{\prime}\in L^2([0,T],\mathbb{C}^*),$ with $\mathbb{C}^*$ denoting the dual space of the space of complex numbers $\mathbb{C}.$}
{The unknown solution $u$ is then approximated} with the finite Fourier series 
\begin{equation}
\label{eq:invFT}
u(t){\approx}\sum\limits_{j=1}^{{N}}\hat{u}_j e^{\imath\omega_j t},\quad t\in[0,T]
\end{equation}
using frequencies
\begin{equation}
\label{eq:freq_wj}
\omega_{j}=2\pi p_j/T,\quad j=1,\dots,N
\end{equation}
from the double-sided spectrum given by $N$-dimensional vector
\begin{equation}
\vec{p}=\big[-\left\lfloor{N/2}\right\rfloor+\delta,\dots,\left\lfloor{N/2}\right\rfloor\big]^{\!\top}\in\mathbb{R}^N,
\end{equation} 
where 
$$\delta=
\begin{cases}
1,& \mathrm{if}\ N\ \mathrm{is\ even}\\ 
0,& \mathrm{if}\ N\ \mathrm{is\ odd}
\end{cases}
$$
and $\lfloor\cdot\rfloor$ denotes the floor function, which returns the greatest integer less than or equal to its argument. 
This yields a finite-dimensional system of equations 
\begin{equation}
\label{eq:var_form_discr}
0=\frac{1}{T}\int\limits_0^T \hat{F}(\hat{\mathbf{u}},t)e^{-\imath\omega_j t}dt,\quad {j=1,\dots,N}
\end{equation}
with respect to the unknown coefficients $\hat{\mathbf{u}}=\left[\hat{u}_{{1}},\dots,\hat{u}_{{N}}\right]^{\!\top}$ in the frequency domain. 
Here $\hat{F}(\hat{\mathbf{u}},t)$ denotes the restriction of $F(u^{\prime},u,t)$ in \eqref{eq:HB_ode} after substitution of $u$ 
from \eqref{eq:invFT} therein. 
Using the calculated frequency components $\hat{u}_j,$ expansion \eqref{eq:invFT} gives the solution in the time domain. 
Note that \IKR{the} approximation of \eqref{eq:FT} using {the left rectangle} quadrature rule on a partition $0=T_0<T_1<\dots<T_{{N-1}}<T$ 
\begin{equation}
\label{eq:discr_FT}
\hat{u}_j=\frac{1}{{N}}\sum\limits_{n=0}^{{N-1}}u(T_n) e^{-\imath\omega_j T_n},\quad {j=1,\dots,N}
\end{equation}
is the discrete Fourier transform (DFT) of vector 
$\mathbf{u}=\left[u(T_0),\dots,u(T_{{N-1}})\right]^{\!\top},$ 
while \eqref{eq:invFT} evaluated at $t=T_n,$ $n=0,\dots,{N-1}$ gives the inverse DFT of $\hat{\mathbf{u}}.$ 
In the following section we exploit the MH discretization idea to solve PP-PC system 
\eqref{eq:PP-PC_sys} for a linear time-periodic problem. 

\section{\IKR{MH solver for linear time-periodic systems}} \label{section:MH_TP}
\IKR{In this section we \IKRrr{describe a combination} of the time-domain solution approaches such as PP-PC \eqref{eq:PP-PC1}-\eqref{eq:PP-PC2} 
and \IKRrr{direct TP discretization} with the frequency-domain representation from Section~\ref{section:Parareal_HB} for linear time-periodic problems, \IKRrr{presented in \cite{Kulchytska-Ruchka_2019ac} and  \cite{Biro_2006aa}, respectively}.}
\subsection{\IKR{Linear PP-PC MH approach}}
\label{sec:lin_pppc_MH}
\IKRrr{Based on \cite{Kulchytska-Ruchka_2019ac} we recall the application of} the MH approach \cite{Biro_2006aa} to parallel-in-time solution of a linear non-autonomous problem of the form
\begin{equation} \label{eq:ode_lin}
\begin{aligned}
  \mathbf{M}\bu^{\prime}(t) + \mathbf{K}\bu(t)&=\bj(t), \quad t \in (0,T),\\
  \bu(0) &= \bu(T)
\end{aligned}
\end{equation}
with matrices $\mathbf{M}$ and $\mathbf{K}$ \IKRrr{having a regular matrix pencil $(\mathbf{M},\mathbf{K})$} 
and RHS $\bj(t)$. 
Discretization by the implicit Euler method on the equidistant coarse grid $0=T_0<T_1<\dots<T_N=T$ defines the 
coarse solution $\mathcal{G}\left(T_n,T_{n-1},\mathbf{U}_{n-1}\right)$ at $T_n$ by 
\begin{equation}
\label{eq:ode_lin_discrete}
\underbrace{\left[\frac{1}{\Delta T}\mathbf{M}+\mathbf{K}\right]}_{=:\mathbf{Q}}\mathcal{G}\left(T_n,T_{n-1},\mathbf{U}_{n-1}\right)=\underbrace{\frac{1}{\Delta T}\mathbf{M}}_{=:\mathbf{C}}\mathbf{U}_{n-1}+\bj(T_n),\quad n=1,\dots,N 
\end{equation}
with step size $\Delta T=T/N.$ 
This yields PP-PC system \eqref{eq:PP-PC_sys} in the explicit matrix-vector form 
\begin{equation}
\label{eq:PP_PC_sys_AVG_lin}
\underbrace{
\begin{bmatrix}
\mathbf{Q} &    &  & {-\mathbf{C}} \\
-\mathbf{C} & \mathbf{Q} &  &  \\
&  \ddots   &  \ddots & \\
&     &  -\mathbf{C}    & \mathbf{Q}
\end{bmatrix}}_{=:\mathbf{G}}
\begin{bmatrix}
	\mathbf{U}_0^{(k+1)}\\
	\mathbf{U}_1^{(k+1)}\\
	\vdots\\
	\mathbf{U}_{N-1}^{(k+1)}
	\end{bmatrix}=
	\underbrace{
		\begin{bmatrix}
\mathbf{r}_N^{(k)}\\
\mathbf{r}_1^{(k)}\\
\vdots\\
\mathbf{r}_{N-1}^{(k)}
\end{bmatrix}}_{=:\mathbf{r}^{(k)}},
\end{equation}
where the RHS is defined by
$$\mathbf{r}_n^{(k)}:=\mathbf{Q}\mathbf{b}_n^{(k)}+\bj(T_n),$$
and vector $\mathbf{b}_n^{(k)},$ $n=1,\dots,N$ is given by \eqref{eq:bnk}.
Originating from \eqref{eq:PP-PC_sys}, system matrix $\mathbf{G}$ has the \IKR{inconvenient} block-cyclic structure. 
The MH solution approach \cite{Biro_2006aa} overcomes this difficulty by decoupling the solution on the coarse level without introducing an additional iteration, in contrast to \eqref{eq:PP-PC_sys_fxd_pt}.
	
Since $\mathbf{U}_0^{(k+1)},\dots,\mathbf{U}_{N-1}^{(k+1)}$ are values of a periodic function at $T_0,\dots,T_{N-1}$ one could 
{express them with the finite Fourier series expansion} \eqref{eq:invFT} discretized at the synchronization points, 
i.e.,  
\begin{equation}
\label{eq:invFT_discr_PPPC}
\mathbf{U}_n^{(k+1)}{=}\sum_{j=1}^N{\hat{\mathbf{U}}_j^{(k+1)}e^{\imath\omega_{j}T_n}},\quad n=0,\dots,N-1,
\end{equation}
with frequencies $\omega_{j}$ given by \eqref{eq:freq_wj}.
Plugging \eqref{eq:invFT_discr_PPPC} into the PP-PC system
\eqref{eq:PP_PC_sys_AVG_lin} and applying DFT \eqref{eq:discr_FT} gives an equivalent system in frequency domain 
%
		\begin{equation}
		\label{eq:PPPC_freq_domain}
			\underbrace{\left(\tilde{\mathbf{F}}\ \mathbf{G}\ \tilde{\mathbf{F}}^{\mathrm{H}}\right)}_{=:\hat{\mathbf{G}}}\hat{\mathbf{U}}^{(k+1)}=\underbrace{\tilde{\mathbf{F}}\mathbf{r}^{(k)}}_{=:\hat{\mathbf{r}}^{(k)}}
		\end{equation}
in terms of the unknown Fourier coefficients $\hat{\mathbf{U}}_{j}^{(k+1)},$ $j=1,\dots,N$ contained in the joint vector $\hat{\mathbf{U}}^{(k+1)}\in{\mathbb{C}}^{Nd}$.  
Matrix $\tilde{\mathbf{F}}$ is defined by $\tilde{\mathbf{F}}=\mathbf{F}\otimes\mathbf{I},$ where $\mathbf{F}$ is the DFT matrix with elements 
\begin{equation}
\label{eq:Fpq}
\mathbf{F}_{jq}=e^{-\imath\omega_{j} T_{q-1}}/{\sqrt{N}},\quad j,q=1,\dots,N,
\end{equation}
$\mathbf{I}$ is a $(d\times d)$-dimensional identity matrix, and `$\otimes$' denotes the Kronecker product of two matrices. $\tilde{\mathbf{F}}^{\mathrm{H}}$ is the Hermite conjugate matrix of $\tilde{\mathbf{F}}$.  
	
We emphasize that {the} transformation from \eqref{eq:PP_PC_sys_AVG_lin} to \eqref{eq:PPPC_freq_domain} comes along with an important property: it transforms the block-cyclic system matrix $\mathbf{G}$ into the block-diagonal system matrix $\hat{\textbf{G}}$. Matrices on the diagonal of $\hat{\textbf{G}}$ can be explicitly calculated as 
\begin{equation}
\label{G_diag_elems}
\hat{\mathbf{G}}_{jj} = \mathbf{Q}-\mathbf{C}e^{-\imath\Delta T\omega_{j}},\quad j=1,\dots,N. 
\end{equation}
\IKRrr{Subsequently, the harmonic components can be found solving}
\begin{equation}
\label{eq:harmonic_indep}
\hat{\mathbf{G}}_{jj}\hat{\mathbf{U}}_j^{(k+1)}=\hat{\mathbf{r}}_j^{(k)},\quad j=1,\dots,N. 
\end{equation}
\IKRrr{Therefore,} additional parallelization \IKRrr{can be obtained on} the coarse level, i.e., solution of $N$ systems of $d$ linear equations could be performed in parallel, without introduction of an additional inner loop. 
Note that matrices $\tilde{\mathbf{F}}$ and $\tilde{\mathbf{F}}^{\mathrm{H}}$ do not have to be explicitly constructed, since {the} system matrices {in \eqref{eq:harmonic_indep}} are determined by \eqref{G_diag_elems} for 
each frequency $\omega_{j}$.

\IKR{Finally, the} solution in the time domain can be obtained by application of the inverse DFT to the calculated solution vector $\hat{\mathbf{U}}^{(k+1)},$ i.e.,
\begin{equation}
\label{eq:FT_TD}
\mathbf{U}^{(k+1)}=\tilde{\mathbf{F}}^{\mathrm{H}}\hat{\mathbf{U}}^{(k+1)}.
\end{equation}
Note that the DFT and its inverse can be efficiently \IKR{applied} using the FFT algorithm \cite{Trefethen_1996aa}, thereby {further} reducing {the} complexity of the transformation.

\subsection{\IKR{Linear TP MH approach}}\label{subsec:TP-FEM_MH_lin}
\IKR{Clearly, the transformation into frequency domain of the form \eqref{eq:PPPC_freq_domain} 
can be also \IKRrr{applied} directly within the discrete TP setting \IKRrr{as it was performed} in \cite{Biro_2006aa}. 
In this case application of implicit Euler's method to \eqref{eq:ode_lin} on the fine grid 
$0=T_0<t_1<\dots<t_{N_{\mathrm{f}}}=T$
\begin{equation}
\label{eq:ode_lin_discrete_TP}
\left[\frac{1}{\delta T}\mathbf{M}+\mathbf{K}\right]\bu_n=\frac{1}{\delta T}\mathbf{M}\bu_{n-1}+\bj(t_n),\quad n=1,\dots,N_{\mathrm{f}}, 
\end{equation}
\IKR{together with the periodicity condition $\bu_0=\bu_{N_{\mathrm{f}}}$} leads to the block-cyclic system
\begin{align}
\label{eq:sys_explicit_lin}
\setlength{\arraycolsep}{2.5pt}
\begin{bmatrix}
\mathbf{C}_{\mathrm{f}}+\mathbf{K} &  &  & -\mathbf{C}_{\mathrm{f}} \\
-\mathbf{C}_{\mathrm{f}} & \mathbf{C}_{\mathrm{f}}+\mathbf{K} &  &  \\
&  \ddots   &  \ddots & \\
&     &  -\mathbf{C}_{\mathrm{f}}    & \mathbf{C}_{\mathrm{f}}+\mathbf{K} 
\end{bmatrix}
\begin{bmatrix}
	\bu_1\\
	\bu_2\\
	\vdots\\
	\bu_{N_{\mathrm{f}}}
	\end{bmatrix}=
		\begin{bmatrix}
\bj(t_1)\\
\bj(t_2)\\
\vdots\\
\bj(t_{N_{\mathrm{f}}})
\end{bmatrix},
\end{align}
where $\mathbf{C}_{\mathrm{f}}:={1}/{\delta T}\cdot\mathbf{M}$ and $\delta T=T/N_{\mathrm{f}}.$ 
Here $\bu_n\approx\bu(t_n)$ denotes the discrete solution at the time instant $t_n,$ $n=1,\dots,N_{\mathrm{f}}.$ 
Solution of system \eqref{eq:sys_explicit_lin} arises in the TP FEM framework, when using the terminology of \cite{Hara_1985aa}. 
We note that the fine 
Euler time step in \eqref{eq:sys_explicit_lin} is much smaller than the coarse one in the PP-PC system 
\eqref{eq:PP_PC_sys_AVG_lin}, i.e., $\delta T\ll \Delta T,$ which means that system 
\eqref{eq:sys_explicit_lin} has a much bigger size compared to that of \eqref{eq:PP_PC_sys_AVG_lin}. }

\IKR{Similarly to \eqref{eq:PPPC_freq_domain} we can apply the MH solver to \eqref{eq:sys_explicit_lin} and compute 
each harmonic coefficient $\hat{\bu}_j,$ $j=1,\dots,N_{\mathrm{f}}$ separately. \IKRrr{This way to solve the time-periodic problem we call TP MH.} In this (linear) case the TP MH \IKRrr{approach} effectively solves 
only one $d$-dimensional linear system, provided sufficient parallelization is possible, 
i.e., when (at least) $N_{\mathrm{f}}$ cores are available. Besides, this approach does not include any 
iteration\IKRrr{,} in contrast to PP-PC \eqref{eq:PP_PC_sys_AVG_lin}, but allows to parallelize the 
computations only through the decoupling of the frequency components.}
\section{\IKR{A novel iterative solver for nonlinear time-periodic systems}}
\label{section:Parareal_new}
\IKR{We now introduce an iterative algorithm which allows to treat also 
nonlinear PP-PC and \IKRrr{discrete} TP formulations in an efficient manner using 
the MH diagonalization technique.}
\subsection{\IKR{Nonlinear PP-PC MH approach}}\label{subsec:nlin_PP-PC_MH}
We now consider the following time-periodic problem for a system of nonlinear 
ODEs
\begin{equation} \label{eq:ode_nlin}
\begin{aligned}
  \mathbf{M}\bu^{\prime}(t) + \mathbf{K}\bigl(\bu(t)\bigr)\bu(t)&=\bj(t), \quad t \in (0,T),\\
  \bu(0) &= \bu(T)
	\end{aligned}
\end{equation}
with nonlinearity present in matrix $\mathbf{K}(\bu(t)),$ 
{e.g., due to nonlinear material characteristics.} 
Similarly to \eqref{eq:PP_PC_sys_AVG_lin}, coarse
discretization with implicit Euler's method on an equidistant grid
leads to the PP-PC system
\begin{align}
\label{eq:PP_PC_sys_explicit_nlin} 
\setlength{\arraycolsep}{2.5pt}
\begin{bmatrix}
\mathbf{Q}_N\bigl(\mathbf{U}_{N-1}^{(k+1)}\bigr) &    &  & {-\mathbf{C}} \\
-\mathbf{C} & \mathbf{Q}_1\bigl(\mathbf{U}_{0}^{(k+1)}\bigr) &  &  \\
&  \ddots   &  \ddots & \\
&     &  -\mathbf{C}    & \mathbf{Q}_{N-1}\bigl(\mathbf{U}_{N-2}^{(k+1)}\bigr) 
\end{bmatrix}
\begin{bmatrix}
	\mathbf{U}_0^{(k+1)}\\
	\mathbf{U}_1^{(k+1)}\\
	\vdots\\
	\mathbf{U}_{N-1}^{(k+1)}
	\end{bmatrix}=
		\begin{bmatrix}
\mathbf{r}_N^{\IKR{k}}\bigl(\mathbf{U}_{N-1}^{(k+1)}\bigr)\\
\mathbf{r}_1^{\IKR{k}}\bigl(\mathbf{U}_{0}^{(k+1)}\bigr)\\
\vdots\\
\mathbf{r}_{N-1}^{\IKR{k}}\bigl(\mathbf{U}_{N-2}^{(k+1)}\bigr)
\end{bmatrix}
\end{align}
with the following definitions
\begin{align*}
\mathbf{Q}_n\bigl(\mathbf{U}_{n-1}^{(k+1)}\bigr)&:=\mathbf{C}+\mathbf{K}\Bigl(\mathcal{G}\bigl(T_n,T_{n-1},\mathbf{U}_{n-1}^{(k+1)}\bigr)\Bigr),\quad \mathbf{C}:={1}/{\Delta T}\cdot\mathbf{M},\\
\mathbf{r}_n^{\IKR{k}}\bigl(\mathbf{U}_{n-1}^{(k+1)}\bigr)&:=\mathbf{Q}_n\bigl(\mathbf{U}_{n-1}^{(k+1)}\bigr)\mathbf{b}_n^{(k)}+\mathbf{j}(T_n),\  n=1,\dots,N,
\end{align*}
where $\mathbf{b}_n^{(k)},$ $n=1,\dots,N$ is given by \eqref{eq:bnk} and $\Delta T=T/N$. 
\IKRrr{We point out the nonlinear dependence of $\mathbf{K}$ on 
$\mathcal{G}\bigl(T_n,T_{n-1},\mathbf{U}_{n-1}\bigr)$ 
in $\mathbf{Q}_n\bigl(\mathbf{U}_{n-1}\bigr)$ in contrast to the implicit Euler discretization, where the nonlinear dependence on the variable $\bu_{n}$ directly is present (see Section~\ref{subsec:TP-FEM_MH_nlin}). This is due to the application of the coarse propagator $\mathcal{G}$ within PP-PC, seen in the operator equation \eqref{eq:PP-PC_sys}.} 
Substituting $\mathcal{G}\bigl(T_n,T_{n-1},\mathbf{U}_{n-1}^{(k+1)}\bigr),$ $n=1,\dots,N$ from PP-PC iteration \eqref{eq:PP-PC1}-\eqref{eq:PP-PC2} 
into \eqref{eq:PP_PC_sys_explicit_nlin} and omitting superscript $k+1$ 
we obtain the following nonlinear system for $\mathbf{U}=\left[\mathbf{U}_0^{\!\top},\dots,\mathbf{U}_{N-1}^{\!\top}\right]^{\!\top}$
\begin{align}
\label{eq:PP_PC_sys_explicit}
\underbrace{
\setlength{\arraycolsep}{2.5pt}
\begin{bmatrix}
\mathbf{Q}\bigl(\mathbf{U}_0-\mathbf{b}_N^{(k)}\bigr) &    &  & {-\mathbf{C}} \\
-\mathbf{C} & \mathbf{Q}\bigl(\mathbf{U}_1-\mathbf{b}_1^{(k)}\bigr) &  &  \\
&  \ddots   &  \ddots & \\
&     &  -\mathbf{C}    & \mathbf{Q}\bigl(\mathbf{U}_{N-1}-\mathbf{b}_{N-1}^{(k)}\bigr) 
\end{bmatrix}
}_{{=:\mathbf{G}^{\IKR{k}}\left(\mathbf{U}\right)}}
\begin{bmatrix}
\mathbf{U}_0\\
\mathbf{U}_1\\
\vdots\\
\mathbf{U}_{N-1}
\end{bmatrix}=
\underbrace{
\begin{bmatrix}
\mathbf{r}_N^{\IKR{k}}\big(\mathbf{U}_{0}\big)\\
\mathbf{r}_1^{\IKR{k}}\big(\mathbf{U}_{1}\big)\\
\vdots\\
\mathbf{r}_{N-1}^{\IKR{k}}\big(\mathbf{U}_{N-1}\big)
\end{bmatrix}
}_{{=:\mathbf{r}^{\IKR{k}}\left(\mathbf{U}\right)}}
\end{align}
where
\begin{align*}
\mathbf{Q}\big(\mathbf{X}\big)&:=\mathbf{C}+\mathbf{K}\left(\mathbf{X}\right),\quad \mathbf{X}\in\mathbb{R}^{d},\quad  \mathbf{C}:={1}/{\Delta T}\cdot\mathbf{M},\\
\mathbf{r}_n^{\IKR{k}}\big(\mathbf{U}_{m}\big)&:=\mathbf{Q}\big(\mathbf{U}_{m}-\mathbf{b}_{n}^{(k)}\big)\mathbf{b}_n^{(k)}+\mathbf{j}(T_n),\quad 1\leq n\leq N,\quad 0\leq m\leq N-1, 
\end{align*}
with $\mathbf{b}_n^{(k)},$ $n=1,\dots,N$ from \eqref{eq:bnk}. 
To solve system \eqref{eq:PP_PC_sys_explicit} at PP-PC iteration $k+1$ we search for the root of mapping $\mathbf{R}^{\IKR{k}}:\mathbb{R}^{Nd}\rightarrow\mathbb{R}^{Nd}$ defined by
\begin{equation}
\label{eq:mappingF}
\mathbf{R}^{\IKR{k}}\left(\mathbf{U}\right)=
\begin{bmatrix}
\mathbf{Q}\bigl(\mathbf{U}_0-\mathbf{b}_N^{(k)}\bigr)\bigl[\mathbf{U}_0-\mathbf{b}_N^{(k)}\bigr]-\mathbf{C}\mathbf{U}_{N-1}-\mathbf{j}(T_N)\\
\mathbf{Q}\bigl(\mathbf{U}_1-\mathbf{b}_1^{(k)}\bigr)\bigl[\mathbf{U}_1-\mathbf{b}_1^{(k)}\bigr] -\mathbf{C}\mathbf{U}_0-\mathbf{j}(T_1) \\
\vdots \\
\mathbf{Q}\bigl(\mathbf{U}_{N-1}-\mathbf{b}_{N-1}^{(k)}\bigr)\bigl[\mathbf{U}_{N-1}-\mathbf{b}_{N-1}^{(k)}\bigr]-\mathbf{C}\mathbf{U}_{N-2}-\mathbf{j}(T_{N-1}) 
\end{bmatrix}
\end{equation}
using the simplified Newton iteration {\cite{Deuflhard_2004aa}}: for $s=0,1,\dots$ \IKR{compute $\mathbf{U}^{(s+1)}=\mathbf{U}^{(k+1,s+1)}$ from} 
\begin{equation}
\label{eq:Newton}
\begin{aligned}
\mathbf{J}_{\mathbf{R}}^{\IKR{k}}\big(\mathbf{U}^{(0)}\big)\mathbf{U}^{(s+1)}=-\Big[\mathbf{R}^{\IKR{k}}\big(\mathbf{U}^{(s)}\big)-\mathbf{J}_{\mathbf{R}}^{\IKR{k}}\big(\mathbf{U}^{(0)}\big)\mathbf{U}^{(s)}\Big],
\end{aligned}
\end{equation}
\IKR{with given initial approximation $\mathbf{U}^{(0)}=\mathbf{U}^{(k+1,0)}.$} 
Jacobian matrix $\mathbf{J}_{\mathbf{R}}^{\IKR{k}}=\displaystyle\frac{\mathrm{d}}{\mathrm{d}\mathbf{U}}\mathbf{R}^{\IKR{k}}$ is defined by 
\begin{equation}
\mathbf{J}_{\mathbf{R}}^{\IKR{k}}\left(\mathbf{U}\right)
=
\setlength{\arraycolsep}{2.5pt}
\begin{bmatrix}
\mathbf{Q}_{\mathrm{d}}\bigl(\mathbf{U}_0-\mathbf{b}_N^{(k)}\bigr) &  &  & -\mathbf{C} \\
-\mathbf{C} & \mathbf{Q}_{\mathrm{d}}\bigl(\mathbf{U}_1-\mathbf{b}_1^{(k)}\bigr) &  &  \\
&  \ddots   &  \ddots & \\
&     &  -\mathbf{C}    & \mathbf{Q}_{\mathrm{d}}\bigl(\mathbf{U}_{N-1}-\mathbf{b}_{N-1}^{(k)}\bigr)
\end{bmatrix} \label{eq:Jacobian_matrix}
\end{equation}
where
\begin{equation}
\label{eq:Qd_Kd}
\mathbf{Q}_{\mathrm{d}}\left(\mathbf{X}\right)=\mathbf{C}+\mathbf{K}_{\mathrm{d}}\left(\mathbf{X}\right),\quad \mathbf{K}_{\mathrm{d}}\left(\mathbf{X}\right)=\frac{\mathrm{d}}{\mathrm{d}\mathbf{X}}\big[\mathbf{K}\left(\mathbf{X}\right)\mathbf{X}\big],\quad\mathbf{X}\in\mathbb{R}^{d}.
\end{equation}
\IKRrr{We note that applying the simplified Newton method we have to construct the system matrix only once at each PP-PC iteration. Now, in order to apply the MH coarse grid correction to the iteration \eqref{eq:Newton} as in Section~\ref{sec:lin_pppc_MH} the block-matrices $\mathbf{Q}_{\mathrm{d}}(\cdot)$ on the diagonal have to be all equal. We therefore choose the} initial approximation for the simplified Newton iteration \eqref{eq:Newton} at PP-PC iteration $k+1$ \IKRrr{as}
\begin{equation}
\label{eq:U0}
\mathbf{U}^{(0)}=\IKR{\mathbf{U}^{(k+1,0)}:=}\Bigl[\bigl(\mathbf{Z}+\mathbf{b}_N^{(k)}\bigr)^{\!\top},\bigl(\mathbf{Z}+\mathbf{b}_1^{(k)}\bigr)^{\!\top}\dots,\bigl(\mathbf{Z}+\mathbf{b}_{N-1}^{(k)}\bigr)^{\!\top}\Bigr]^{\!\top}
\end{equation}
for a given vector $\mathbf{Z}\in\mathbb{R}^d.$ \IKRrr{The}  
simplified Newton iteration \eqref{eq:Newton} \IKRrr{then} reads
\begin{equation}
\label{eq:Newton_init0}
\mathbf{G}_{\mathrm{d}}\mathbf{U}^{(s+1)}=\mathbf{h}^{(s)} 
\end{equation}
with system matrix
\begin{equation}
\label{eq:Newton_init0_matr}
\mathbf{G}_{\mathrm{d}}\IKR{:=\mathbf{J}_{\mathbf{R}}^{k}\bigl(\mathbf{U}^{(0)}\bigr)}=
\setlength{\arraycolsep}{2.5pt}
\begin{bmatrix}
\mathbf{Q}_{\mathrm{d}}\left(\mathbf{Z}\right) &  &  & -\mathbf{C} \\
-\mathbf{C} & \mathbf{Q}_{\mathrm{d}}\left(\mathbf{Z}\right) &  &  \\
&  \ddots   &  \ddots & \\
&     &  -\mathbf{C}    & \mathbf{Q}_{\mathrm{d}}\left(\mathbf{Z}\right)
\end{bmatrix}
\end{equation}
and RHS 
\begin{equation}
\label{eq:hs}
\mathbf{h}^{(s)}:=-\bigl[\mathbf{R}^{\IKR{k}}\bigl(\mathbf{U}^{(s)}\bigr)-\mathbf{G}_{\mathrm{d}}\mathbf{U}^{(s)}\bigr].
\end{equation}
The matrix \IKRrr{in \eqref{eq:Newton_init0_matr}} remains constant over the Newton iterations and has the same 
block-cyclic structure as that of the PP-PC system \eqref{eq:PP_PC_sys_AVG_lin} for a linear problem. 
Hence, a MH solver {can} be also applied to the Newton system \eqref{eq:Newton_init0}. 
Analogously to \eqref{eq:PPPC_freq_domain} we obtain
\begin{equation}
\label{eq:PPPC_Newton_freq_domain}
\underbrace{\left(\tilde{\mathbf{F}}\ \mathbf{G}_{\mathrm{d}} \ \tilde{\mathbf{F}}^{\mathrm{H}}\right)}_{=:\hat{\mathbf{G}}_{\mathrm{d}} }\hat{\mathbf{U}}^{(s+1)}
=\tilde{\mathbf{F}}\mathbf{h}^{(s)},  
\end{equation}
where matrix $\hat{\mathbf{G}}_{\mathrm{d}}$ is block-diagonal with each $(d\times d)$-dimensional block given by
\begin{equation}
\label{G_diag_elems_Newton}
\big[\hat{\mathbf{G}}_{\mathrm{d}}\big]_{jj} = \mathbf{Q}_{\mathrm{d}}\left(\mathbf{Z}\right)-\mathbf{C}e^{-\imath\Delta T\omega_{j}},\quad j=1,\dots,N
\end{equation}
and {the} DFT matrix $\tilde{\mathbf{F}}$ is as previously defined by \eqref{eq:Fpq}. 
{The} frequency domain solution $\hat{\mathbf{U}}^{(s+1)}$ at Newton iteration $s+1$ is then transformed into the time domain by the inverse DFT, i.e.,
\begin{equation}
\label{eq:FT_TD_Newton}
\mathbf{U}^{(s+1)}=\tilde{\mathbf{F}}^{\mathrm{H}}\hat{\mathbf{U}}^{(s+1)}.
\end{equation}
\begin{rem}
	We note that {the} simplified Newton iteration \eqref{eq:Newton} could be also interpreted as a linear iterative method based on an additive splitting of the system matrix in \eqref{eq:PP_PC_sys_explicit}. Indeed, using {any} constant matrix {$\mathbf{H}^{k}\in\mathbb{R}^{Nd\times Nd}$} we can introduce \IKR{at PP-PC iteration $k+1$} a linear iteration for $s=0,1,\dots$ 
	\begin{equation}
	\label{eq:add_splitting}
	{\mathbf{H}}^{k}\mathbf{U}^{(s+1)}=\left[{\mathbf{H}}^{k}-\mathbf{G}^{\IKR{k}}\big(\mathbf{U}^{(s)}\big)\right]\mathbf{U}^{(s)}+\mathbf{r}^{\IKR{k}}\big(\mathbf{U}^{(s)}\big),
	\end{equation}
	which in case of $\mathbf{H}^{k}\!:=\mathbf{J}_{\mathbf{R}}^{\IKR{k}}\big(\mathbf{U}^{(0)}\big)$ becomes exactly the proposed simplified Newton iteration \eqref{eq:Newton}. \IKRrr{We note that that the choice 
	$$\mathbf{H}^{k}(\mathbf{U})=\mathrm{diag}\Bigl[\mathbf{Q}_{\mathrm{d}}\bigl(\mathbf{U}_0-\mathbf{b}_N^{(k)}\bigr),\mathbf{Q}_{\mathrm{d}}\bigl(\mathbf{U}_1-\mathbf{b}_1^{(k)}\bigr),\dots,\mathbf{Q}_{\mathrm{d}}\bigl(\mathbf{U}_{N-1}-\mathbf{b}_{N-1}^{(k)}\bigr)\Bigr]$$
	yields Gander's Jacobi-like fixed-point iteration \eqref{eq:PP-PC_sys_fxd_pt}.} The idea of choosing a modified fixed point iteration operator to ensure a convenient block structure is well known and is for example applied by Biro \cite{Biro_2006aa} and Gander \cite{Gander_2017aa}. However, both do not discuss convergence.	
\end{rem}
\IKRrr{
\begin{rem}
We point out that the MH approach is valid not only when applying the implicit Euler method but also other time intergation algorithms. 
For instance, using the implicit trapezoidal rule (also called the Crank-Nicolson method) in the linear case, the system \eqref{eq:PP_PC_sys_AVG_lin} maintains the same block-cyclic structure with the following definitions of the block matrices and the RHS:
$$\mathbf{C}:=\frac{1}{\Delta T}\mathbf{M}-0.5\mathbf{K},\quad\mathbf{Q}:=\mathbf{C}+\mathbf{K},\quad\mathbf{r}_n^{(k)}:=\mathbf{Q}\mathbf{b}_n^{(k)}+0.5\left[\bj(T_n)+\bj(T_{n-1})\right].$$
\end{rem}
}

\subsection{\IKR{Nonlinear TP MH approach}}
\label{subsec:TP-FEM_MH_nlin}
Similarly as in Section~\ref{subsec:TP-FEM_MH_lin}, one could exploit the proposed methodology 
to solve a time-periodic problem, without applying a parallel-in-time method. 
More specifically, the fine discretization of the periodic system \eqref{eq:ode_nlin} on $0=T_0<t_1<\dots<t_{N_{\mathrm{f}}}=T$
using implicit Euler's method 
with time step size $\delta T=T/N_{\mathrm{f}}$ leads to the system of nonlinear algebraic equations
\begin{align}
\label{eq:sys_explicit_nlin}
\setlength{\arraycolsep}{2.5pt}
\begin{bmatrix}
\mathbf{C}_{\mathrm{f}}+\mathbf{K}\left(\bu_1\right) &  &  & -\mathbf{C}_{\mathrm{f}} \\
-\mathbf{C}_{\mathrm{f}} & \mathbf{C}_{\mathrm{f}}+\mathbf{K}\left(\bu_2\right) &  &  \\
&  \ddots   &  \ddots & \\
&     &  -\mathbf{C}_{\mathrm{f}}  & \mathbf{C}_{\mathrm{f}}+\mathbf{K}\left(\bu_{N_{\mathrm{f}}}\right) 
\end{bmatrix}
\begin{bmatrix}
	\bu_1\\
	\bu_2\\
	\vdots\\
	\bu_{N_{\mathrm{f}}}
	\end{bmatrix}=
		\begin{bmatrix}
\bj(t_1)\\
\bj(t_2)\\
\vdots\\
\bj(t_{N_{\mathrm{f}}})
\end{bmatrix},
\end{align}
where $\mathbf{C}_{\mathrm{f}}:={1}/{\delta T}\cdot\mathbf{M}.$ 
Choosing vector $\mathbf{z}\in\mathbb{R}^d$ in the initial approximation 
\begin{equation}
\label{eq:U0_TP}
\bu^{(0)}:=\Bigl[\bz^{\!\top},\bz^{\!\top}\dots,\bz^{\!\top}\Bigr]^{\!\top}{\in\mathbb{R}^{N_{\mathrm{f}}d}}
\end{equation}
we obtain the simplified Newton iteration of form \eqref{eq:Newton_init0} applied to \eqref{eq:sys_explicit_nlin}
\begin{equation}
\label{eq:Newton_init0_TP}
\underbrace{
\setlength{\arraycolsep}{2.5pt}
\begin{bmatrix}
\mathbf{Q}_{\mathrm{f}}\left(\bz\right) &  &  & -\mathbf{C}_{\mathrm{f}} \\
-\mathbf{C}_{\mathrm{f}} & \mathbf{Q}_{\mathrm{f}}\left(\bz\right) &  &  \\
&  \ddots   &  \ddots & \\
&     &  -\mathbf{C}_{\mathrm{f}}  & \mathbf{Q}_{\mathrm{f}}\left(\bz\right)
\end{bmatrix}
}_{=:\mathbf{G}_{\mathrm{d}}}
\begin{bmatrix}
\bu_1^{(s+1)}\\
\bu_2^{(s+1)}\\
\vdots\\
\bu_{N_{\mathrm{f}}}^{(s+1)}
\end{bmatrix}
=
\underbrace{
\begin{bmatrix}
\bigl[\mathbf{K}_{\mathrm{d}}\left(\bz\right)-\mathbf{K}\bigl(\bu_1^{(s)}\bigr)\bigr]\bu_1^{(s)}+\bj_1\\
\bigl[\mathbf{K}_{\mathrm{d}}\left(\bz\right)-\mathbf{K}\bigl(\bu_2^{(s)}\bigr)\bigr]\bu_2^{(s)}+\bj_2\\
\vdots\\
\bigl[\mathbf{K}_{\mathrm{d}}\left(\bz\right)-\mathbf{K}\bigl(\bu_{N_{\mathrm{f}}}^{(s)}\bigr)\bigr]\bu_{N_{\mathrm{f}}}^{(s)}+\bj_{N_{\mathrm{f}}}
\end{bmatrix}
}_{=:\mathbf{h}^{(s)}},
\end{equation}
with $\mathbf{Q}_{\mathrm{f}}\left(\bz\right):=\mathbf{C}_{\mathrm{f}}+\mathbf{K}_{\mathrm{d}}\left(\bz\right)$ 
and $\mathbf{K}_{\mathrm{d}}$ defined in \eqref{eq:Qd_Kd}. 
One could apply \IKR{the} MH solver \eqref{eq:PPPC_Newton_freq_domain} 
at Newton iteration $s+1$ and solve a linear algebraic system for 
each harmonic coefficient $\hat{\bu}_j^{(s+1)},$ $j=1,\dots,N_{\mathrm{f}}$ 
separately. \todo{This is to our best knowledge essentially the same idea 
as the fixed point method proposed in \cite{Biro_2006aa}.}

{When many central processing units (CPUs) are available} one may expect the parallelization of 
the {MH} frequency domain solution at each Newton iteration \eqref{eq:Newton_init0_TP} {to}  
outperform the speed up provided by the parallel-in-time solution on the fine grid {in PP-PC \eqref{eq:PP_PC_sys_explicit_nlin}}. 
However, we note that \IKR{the} application of the MH solver to the PP-PC system and not to {\eqref{eq:sys_explicit_nlin}} 
is especially beneficial when the number of CPUs is limited ($<N_{\mathrm{f}}$), that is, 
when one could not calculate each harmonic coefficient $\hat{\bu}_j^{(s+1)},$ $j=1,\dots,N_{\mathrm{f}}$ 
in parallel.

\subsection{Pseudocode and implementation details}
\label{subsec:pseudocode}
The error used in the termination criterion 
\IKR{for the periodic Parareal-based algorithms (PP-IC and PP-PC) calculates the maximum among the jumps at 
the synchronization points $T_n,$ with $n=1,\dots,N-1,$ as well as the periodicity jump of the solution at $T_0$ and 
$T_N$, i.e., the error at the periodic Parareal (PP) iteration $k$ is computed as
\begin{align}
\label{eq:error_tol_PP}
\varepsilon_{\mathrm{PP}}^{(k)}=\max\left\{\max\limits_{1\leq n\leq N-1}\|\mathbf{U}^{(k)}_n-\mathcal{F}\bigl(T_n,T_{n-1},\mathbf{U}^{(k)}_{n-1}\bigr)\|_{*},\ 
\|\mathbf{U}^{(k)}_0-\mathcal{F}\bigl(T_N,T_{N-1},\mathbf{U}^{(k)}_{N-1}\bigr)\|_{*}\right\}, 
\end{align}
where for two $d$-dimensional vectors $\mathbf{u}$ and $\mathbf{v}$ the $\|\cdot\|_{*}$-norm is defined by
\begin{equation}
\label{eq:error_tol}
\|\mathbf{u}-\mathbf{v}\|_{*}=\frac{\|\mathbf{u}-\mathbf{v}\|_{2}}{\text{aTol}+\text{rTol}\|\mathbf{u}\|_{2}},
\end{equation}
and $\text{aTol}=10^{-6}$ and $\text{rTol}=10^{-3}$ are chosen as the absolute and the relative tolerances, respectively. 
Here $\|\cdot\|_{2}$ denotes the Euclidean norm in the $d$-dimensinal vector space.} 

\IKR{Furthermore, at PP iteration $k+1$ the errors of the inner fixed point \eqref{eq:PP-PC_sys_fxd_pt} 
and the simplified Newton \eqref{eq:Newton} iterations 
are obtained 
from the difference of the coarse solution vectors at the two subsequent iterations $s+1$ and $s$ for $s=0,1,\dots$ through 
\begin{equation}
\label{eq:error_tol_IT}
\varepsilon_{\mathrm{PP,\ it}}^{(k+1,s+1)}=\max\limits_{0\leq n\leq N-1}\|\mathbf{U}^{(k+1,s+1)}_n-\mathbf{U}^{(k+1,s)}_{n}\|_{*}.
\end{equation}
This is similar tothe error of the simplified Newton iteration \eqref{eq:Newton_init0_TP}, 
applied directly to the discrete TP formulation, where the error 
is defined from the (fine) solutions by
\begin{equation}
\label{eq:error_tol_IT_TP}
\varepsilon_{\mathrm{TP,\ it}}^{(s+1)}=\max\limits_{1\leq n\leq N_{\mathrm{f}}}\|\bu^{(s+1)}_n-\bu^{(s)}_{n}\|_{*}.
\end{equation}
Finally, for $k\geq 1$ the classical sequential time-stepping considers the periodicity error 
\begin{equation}
\label{eq:error_tol_seq}
\varepsilon_{\mathrm{seq}}(k)=\|\bu(T_0+kT)-\bu(T_0+(k-1)T)\|_{*}
\end{equation}
after each period $T$ until the periodic steady state is reached. 
The numerical calculations 
are terminated when the values of the corresponding errors from 
\eqref{eq:error_tol_PP}, \eqref{eq:error_tol_IT}-\eqref{eq:error_tol_seq} become smaller than $1$.} 

\IKRrr{The considered methods are implemented in GNU Octave version 3.8.2 using the `parallel' package\footnotemark[1]\footnotetext[1]{\url{http://octave.sourceforge.io/parallel}} for parallelization. We use Octave's Fourier analysis implemented by the parallel version of FFTW3\footnotemark[2]\footnotetext[2]{\url{http://www.fftw.org/parallel}} using $80$ threads. The code is executed on an
Intel Xeon cluster with $\SI{80\times2.00}{\giga\hertz}$ cores, i.e., $8\times$E7-8850 and 1TB DDR3 memory.}

\algblock{ParFor}{EndParFor}
\algnewcommand\algorithmicparfor{\textbf{parfor}}
\algnewcommand\algorithmicpardo{\textbf{do}}
\algnewcommand\algorithmicendparfor{\textbf{end\ parfor}}
\algrenewtext{ParFor}[1]{\algorithmicparfor\ #1\ \algorithmicpardo}
\algrenewtext{EndParFor}{\algorithmicendparfor}
\begin{algorithm}[t]
	\caption{Proposed PP-PC MH algorithm from Section~\ref{subsec:nlin_PP-PC_MH}}
	\label{alg:pppc_mh}
	\begin{algorithmic}[1]
		\State initialize: $\overline{\mathbf{u}}^{(0)}, \tilde{\mathbf{u}}^{(0)}\leftarrow \mathbf{0}$, 
		$\varepsilon_{\mathrm{PP}}^{(0)}\leftarrow 1,$ set counter: $k\leftarrow 0$;
		\While{$k \le {K}\, {\mathbf{and}}\; \varepsilon_{\mathrm{PP}}^{(k)}\geq 1$}		
		\State calculate: $\mathbf{b}^{(k)}_j\leftarrow \tilde{\mathbf{u}}_j^{(k)}-\overline{\mathbf{u}}_j^{(k)},$ choose $\mathbf{Z}$, assign $\hat{\mathbf{G}}_{\mathrm{d}}$ by \eqref{G_diag_elems_Newton};
		\State initialize: $\mathbf{U}^{(k+1,0)}$ by \eqref{eq:U0}, $\varepsilon_{\mathrm{PP,\ it}}^{(k+1,0)}\leftarrow 1,$ set counter: $s\leftarrow 0$;
		\While{$s \le {S}\, {\mathbf{and}}\; \varepsilon_{\mathrm{PP,\ it}}^{(k+1,s)}\geq 1$}	
		\State {update the RHS $\mathbf{h}^{(s)}$ by \eqref{eq:hs};}
		\ParFor{$j \leftarrow 1 ,\, N$}
		\State FFT of the RHS: $\hat{\mathbf{h}}^{(s)}_j\leftarrow\left(\left[\mathbf{F}_{j1}\dots\mathbf{F}_{jN}\right]\otimes\mathbf{I}\right)\mathbf{h}^{(s)}$; 
		\State solve: $[\hat{\mathbf{G}}_{\mathrm{d}}\big]_{jj}\hat{\mathbf{U}}^{(k+1,s+1)}_j=\hat{\mathbf{h}}^{(s)}_j$;
		\State inverse FFT: $\mathbf{U}^{(k+1,s+1)}_j\leftarrow\left(\left[\mathbf{F}^{\mathrm{H}}_{j1}\dots \mathbf{F}^{\mathrm{H}}_{jN}\right]\otimes\mathbf{I}\right)\hat{\mathbf{U}}^{(k+1,s+1)}_j$;
		\EndParFor
		\State update $\varepsilon_{\mathrm{PP,\ it}}^{(k+1,s+1)}$ by \eqref{eq:error_tol_IT}, increment counter: $s \leftarrow s+1$;
		\EndWhile
		\ParFor{$j \leftarrow 1 ,\, N$}
		\State solve fine problem: $\tilde{\mathbf{u}}_j^{(k+1)}\leftarrow \mathcal{F}(T_j, T_{j-1}, \mathbf{U}_{j-1}^{(k+1,s)})$;
		\State solve coarse problem: $\overline{\mathbf{u}}_j^{(k+1)}\leftarrow \mathcal{G}(T_j, T_{j-1}, \mathbf{U}_{j-1}^{(k+1,s)})$;
		\EndParFor
		\State update $\varepsilon_{\mathrm{PP}}^{(k+1)}$ by \eqref{eq:error_tol_PP}, increment counter: $k \leftarrow k+1$;
		\EndWhile
	\end{algorithmic}
\end{algorithm}	
\IKRrr{In Algorithm~\ref{alg:pppc_mh} we present the pseudocode of the nonlinear PP-PC MH method, proposed in Section~\ref{subsec:nlin_PP-PC_MH}.} 

\subsection{Convergence analysis \IKR{of the proposed algorithm}} 
Convergence of the simplified Newton method \eqref{eq:Newton} is determined by the result presented in \cite[Theorem 2.5]{Deuflhard_2004aa}, which originates from \cite{Ortega_2000aa}.
We recall the theorem here in terms of the notations introduced above.
\begin{thm}[Convergence of the simplified Newton method \cite{Deuflhard_2004aa}] \label{thm:Newton}
	Let $\mathbf{R}:D\rightarrow\mathbb{R}^{N}$ be a continuously differentiable mapping with $D\subset\mathbb{R}^{N}$ open and convex. 
	Let $\mathbf{U}^{(0)}\in D$ denote a given starting point so that $\mathbf{J}_{\mathbf{R}}^{\IKR{k}}\left(\mathbf{U}^{(0)}\right)$ is invertible. Assume 
	the affine covariant Lipschitz condition
	\begin{equation}
	\label{eq:cond_Newton}
	\left\|\mathbf{J}_{\mathbf{R}}^{\IKR{k}}\bigl(\mathbf{U}^{(0)}\bigr)^{-1}
	\left[\mathbf{J}_{\mathbf{R}}^{\IKR{k}}\left(\mathbf{U}\right)-\mathbf{J}_{\mathbf{R}}^{\IKR{k}}\bigl(\mathbf{U}^{(0)}\bigr)\right]\right\|
	\leq\delta_0\bigl\|\mathbf{U}-\mathbf{U}^{(0)}\bigr\|
	\end{equation}
	holds for all $\mathbf{U}\in D.$ Let 
	\begin{equation}
	\label{eq:cond_Newton_12}
	h_0:=\delta_0\left\|\mathbf{U}^{(1)}-\mathbf{U}^{(0)}\right\|\leq0.5
	\end{equation}
	and assume that the closure $\bar{S}\left(\mathbf{U}^{(0)},\rho\right)$ of a ball with center in $\mathbf{U}^{(0)}$ and radius $\rho$  
	is a subset of $D,$ with 
	\begin{equation}
	\label{eq:rho_Newton}
	\rho=\left(1-\sqrt{1-2h_0}\right)/\delta_0.
	\end{equation}
	Then the simplified Newton iterates $\mathbf{U}^{(s+1)},$ $s=0,1,\dots$ generated by \eqref{eq:Newton} remain in $\bar{S}\left(\mathbf{U}^{(0)},\rho\right)$ and converge to some $\mathbf{U}^*$ with $\mathbf{R}^{\IKR{k}}\left(\mathbf{U}^*\right)=\mathbf{0}.$ 
\end{thm}

\subsubsection{Numerical analysis of a model problem} 
We {discuss the} convergence {analysis of the} simplified Newton iteration \eqref{eq:Newton} 
for the PP-PC algorithm applied to a time-periodic 
problem for a single ($d=1$) {ODE}. \IKRrr{Motivated by the properties of the nonlinear material curves proposed in \cite{Heise_1994aa} we consider the following model problem} 
\begin{equation} \label{eq:ode_1d}
\begin{aligned}
mu^{\prime}(t) + \kappa(|u|)u(t)&=j(t), \quad t \in (0,T),\\
u(0) &= u(T),
\end{aligned}
\end{equation}
with $m\in\mathbb{R}_{0}^{+},$ $\kappa:\mathbb{R}_{0}^{+}\rightarrow\mathbb{R}_{0}^{+},$ 
unknown function $u:[0,T]\rightarrow\mathbb{R},$ and $T$-periodic input $j,$ i.e., $j(0)=j(T).$ 
We assume that the nonlinear function $\kappa$ satisfies the following conditions.
\begin{ass}
\label{thm:kappa_ass}
Let for $\kappa:\mathbb{R}_0^+\rightarrow\mathbb{R}_0^+$ it holds:
\begin{itemize}
	\item[$\bullet$]  $\kappa(t)\geq c_1>0,$ $t\in\mathbb{R}_0^+,$
	\item[$\bullet$]  function $\IKRrr{f}:t\rightarrow\kappa(t)t$ is strongly monotone with monotonicity constant $c_1,$ i.e., 
	\begin{equation}
	\label{eq:str_monotone}
	(\IKRrr{f(t)-f(s)})(t-s)\geq c_1(t-s)^2\quad \forall t,s\in\mathbb{R}_0^+,
	\end{equation}
	\item[$\bullet$] $\kappa\in C^1(\mathbb{R}_0^+)$ and $\lim\limits_{t\to\infty}\kappa^{\prime}(t)=0,$  
	\item[$\bullet$] function $\IKRrr{g}:t\rightarrow\kappa^{\prime}(t)t$ is Lipschitz continuous, i.e.,
	\begin{equation}
	\label{eq:Lipschitz_cont}
	\exists L_1>0:\quad |\IKRrr{g(t)-g(s)}|\leq L_1|t-s|\quad \forall t,s\in\mathbb{R}_0^+.
	\end{equation}
\end{itemize} 
\end{ass}
We now apply the simplified Newton iteration \eqref{eq:Newton} to \eqref{eq:ode_1d}. Analogously to \eqref{eq:Jacobian_matrix} the Jacobian matrix for the one-dimensional (1D) problem \eqref{eq:ode_1d} is determined at PP-PC iteration $k$ by
\begin{equation}
\mathbf{J}_{\mathbf{R}}^{\IKR{k}}\left(\mathbf{U}\right)=
\setlength{\arraycolsep}{2.5pt}
\begin{bmatrix}
c+\kappa_{\mathrm{d}}\bigl(u_0-b_N^{(k)}\bigr) &  &  & -c \\
-c & c+\kappa_{\mathrm{d}}\bigl(u_1-b_1^{(k)}\bigr) &  &  \\
&  \ddots   &  \ddots & \\
&     &  -c    & c+\kappa_{\mathrm{d}}\bigl(u_{N-1}-b_{N-1}^{(k)}\bigr)
\end{bmatrix},
\end{equation}
where $\mathbf{U}=\left[u_0,\dots,u_{N-1}\right]^{\!\top},$ $c=m/\Delta T$ {with $\Delta T=T/N,$} and function $\kappa_{\mathrm{d}}:\mathbb{R}\mapsto\mathbb{R}_0^+$ given by 
\begin{equation}
\label{eq:kappad}
\kappa_{\mathrm{d}}(x)=
\kappa^{\prime}(|x|)|x|+\kappa(|x|),\ x\in\mathbb{R}. 
\end{equation}
As in \eqref{eq:U0} choosing a fixed value $z\in\mathbb{R}$ we define
\begin{equation}
\label{eq:Ub}
\mathbf{U}^{(0)}=\IKR{\mathbf{U}^{(k+1,0)}:=}\left[z+b_N^{(k)},z+b_1^{(k)}\dots,z+b_{N-1}^{(k)}\right]^{\!\top}
\end{equation}
to be the initial approximation for the Newton algorithm at PP-PC iteration $k+1$. 
{We denote the} Jacobian matrix {$\mathbf{J}_{\mathbf{R}}^{\IKR{k}}$} 
at the chosen $\mathbf{U}^{(0)}$ {from \eqref{eq:Ub} by $\mathbf{G}_{\mathrm{d}}$  
which is given by}
\begin{equation}
\label{eq:Gd_JacU0}
\mathbf{G}_{\mathrm{d}}=
\setlength{\arraycolsep}{2.5pt}
\begin{bmatrix}
c+\kappa_{\mathrm{d}}\left(z\right) &  &  & -c \\
-c & c+\kappa_{\mathrm{d}}\left(z\right) &  &  \\
&  \ddots   &  \ddots & \\
&     &  -c     & c+\kappa_{\mathrm{d}}\left(z\right)
\end{bmatrix}.
\end{equation}
Before demonstrating the convergence results based on the Theorem~\ref{thm:Newton} 
we derive some properties of the function $\kappa_{\mathrm{d}}$ 
in the following lemma.
\begin{lem}
\label{thm:kappa_bdd_Lipschitz}
Under the Assumption~\ref{thm:kappa_ass} on $\kappa$ we have that function $\kappa_{\mathrm{d}}:\mathbb{R}\mapsto\mathbb{R}_0^+$ defined in \eqref{eq:kappad} is
\begin{itemize}
	\item[1)] 
	bounded from below by $c_1$, i.e.,
	\begin{equation}
	\label{eq:coercive}
	\kappa_{\mathrm{d}}(x)\geq c_1\quad \forall x\in\mathbb{R},
	\end{equation}
	\item[2)]
	Lipschitz continuous with Lipschitz constant $L_2=\kappa_2+L_1,$ 
	\IKRrr{where $\kappa_2:=\sup\limits_{t\in\mathbb{R}_0^+}|\kappa^{\prime}(t)|,$} i.e.,
	\begin{equation}
	\label{eq:Lipschitz_kappa_d}
	|\kappa_{\mathrm{d}}(x)-\kappa_{\mathrm{d}}(y)|\leq L_2|x-y|\quad \forall x,y\in\mathbb{R},
	\end{equation}
\end{itemize}
\begin{proof} The proof is based on \cite{Heise_1994aa}, \cite{Pechstein_2004aa}.
\begin{itemize}
	\item[1)]  
	\IKRrr{Based on the proof of Lemma~2.8 in \cite{Pechstein_2004aa} strong monotonicity of $f:t\rightarrow\kappa(t)t,$ $t\in\mathbb{R}_0^+$ with constant $c_1$ implies strong monotonicity of $f_{\text{abs}}:t\rightarrow\kappa(|t|)t,$ $t\in\mathbb{R}$ with the same constant $c_1,$ which yields $\kappa_{\mathrm{d}}(t)=f_{\text{abs}}^{\prime}(t)\geq c_1,$ $t\in\mathbb{R}.$}
	\item[2)] 
	\IKRrr{First, continuity of $\kappa^{\prime}$ on every closed interval $[0,q],$ $0<q<\infty$ and its vanishing on infinity implies its boundedness, i.e., $|\kappa^{\prime}(t)|\leq\kappa_2<\infty,$ $t\in\mathbb{R}_0^+.$ Then using the mean value theorem we have that $\kappa$ is Lipschitz continuous with Lipschitz constant $\kappa_2,$ i.e.,
	\begin{equation}
	\label{eq:Lipschitz_kappa}
	|\kappa(t)-\kappa(s)|=|\kappa^{\prime}(\xi)(t-s)|\leq \kappa_2|t-s|\quad \mathrm{for\ some}\ \xi\in\mathbb{R}_0^+\ \textrm{and}\ \forall t,s\in\mathbb{R}_0^+. 
	\end{equation}
	Using this and the assumption on $g$, we see that $\kappa_{\mathrm{d}}$ is Lipschitz continuous with Lipschitz constant $L_2=\kappa_2+L_1.$ 
	}
	\end{itemize}
\end{proof}
\end{lem}

\begin{thm}
\label{thm:conv_kappad}
If the nonlinear function $\kappa$ in \eqref{eq:ode_1d} satisfies the Assumption~\ref{thm:kappa_ass},
then {the} affine covariant Lipschitz condition \eqref{eq:cond_Newton} holds 
with constant $\delta_0={L_2}/{c_1},$ where $L_2$ is the Lipschitz constant of $\kappa_{\mathrm{d}}.$
%
\begin{proof}
First, we observe that for a chosen $\mathbf{U}^{(0)}$ in \eqref{eq:Ub} the Jacobian $\mathbf{J}_{\mathbf{R}}^{\IKR{k}}\bigl(\mathbf{U}^{(0)}\bigr)=\mathbf{G}_{\mathrm{d}}$ given in \eqref{eq:Gd_JacU0} can be decomposed as 
\begin{equation}
\label{eq:svd}
\mathbf{G}_{\mathrm{d}}={\mathbf{F}}^{\mathrm{H}}\hat{\mathbf{G}}_{\mathrm{d}}{\mathbf{F}},
\end{equation}
with unitary DFT matrix $\mathbf{F}$ defined in \eqref{eq:Fpq} 
and diagonal matrix $\hat{\mathbf{G}}_{\mathrm{d}}$ \IKRrr{defined in \eqref{G_diag_elems_Newton}} containing eigenvalues of $\mathbf{G}_{\mathrm{d}}$. 
This implies
\begin{equation}
\label{eq:norm_Gd}
\left\|\mathbf{G}_{\mathrm{d}}^{-1}\right\|_2=
\left\|\left({\mathbf{F}}^{\mathrm{H}}\hat{\mathbf{G}}_{\mathrm{d}}{\mathbf{F}}\right)^{-1}\right\|_2=
\left\|\hat{\mathbf{G}}_{\mathrm{d}}^{-1}\right\|_2=
\max\limits_{1\leq j\leq N}\left|\Bigl[\kappa_{\mathrm{d}}\left(z\right)+\bigl(1-e^{-\imath\Delta T\omega_{j}}\bigr)c\Bigr]^{-1}\right|,
\end{equation}
where $\|\cdot\|_2=\sigma_{\mathrm{max}}(\cdot)$ denotes the largest singular value of a matrix, 
which for a normal matrix such as $\mathbf{G}_{\mathrm{d}}^{-1}$ is equal to its spectral radius. 
Using \eqref{eq:norm_Gd} we have {for the Lipschitz condition \eqref{eq:cond_Newton} that}
\begin{equation}
\label{eq:estim}
\begin{aligned}
&\left\|\mathbf{J}_{\mathbf{R}}^{\IKR{k}}\bigl(\mathbf{U}^{(0)}\bigr)^{-1}\left[\mathbf{J}_{\mathbf{R}}^{\IKR{k}}\bigl(\mathbf{U}\bigr)-\mathbf{J}_{\mathbf{R}}^{\IKR{k}}\bigl(\mathbf{U}^{(0)}\bigr)\right]\right\|_2
\leq
\left\|\mathbf{G}_{\mathrm{d}}^{-1}\right\|_2\left\|\mathbf{J}_{\mathbf{R}}^{\IKR{k}}\bigl(\mathbf{U}\bigr)-\mathbf{J}_{\mathbf{R}}^{\IKR{k}}\bigl(\mathbf{U}^{(0)}\bigr)\right\|_2=\\
&
\max\limits_{1\leq j\leq N}\left|\Bigl[\kappa_{\mathrm{d}}\left(z\right)+\bigl(1-e^{-\imath\Delta T\omega_{j}}\bigr)c\Bigr]^{-1}\right|
\max\limits_{\substack{1\leq n\leq N \\ 0\leq m\leq N-1}}\left|\kappa_{\mathrm{d}}\bigl(u_m-b_n^{(k)}\bigr)-\kappa_{\mathrm{d}}\left(z\right)\right|.
\end{aligned}
\end{equation}
Now, to show convergence of the simplified Newton method \eqref{eq:Newton} one has to prove that expression \eqref{eq:estim} can be estimated by
$\delta_0\left\|\mathbf{U}-\mathbf{U}^{(0)}\right\|_2,$ $\delta_0>0$ from above.
Indeed, 
since due to the Lemma~\ref{thm:kappa_bdd_Lipschitz} boundedness condition \eqref{eq:coercive} holds for $\kappa_{\mathrm{d}}$
and since $c\in\mathbb{R}_{0}^{+},$ we have for $1\leq j\leq N$
\begin{equation}
\label{eq:bdd_inv}
\begin{aligned}
&\left|\Bigl[\kappa_{\mathrm{d}}\left(z\right)+\bigl(1-e^{-\imath\Delta T\omega_{j}}\bigr)c\Bigr]^{-1}\right|\leq 
\Bigl[\Re\Bigl(\kappa_{\mathrm{d}}\left(z\right)+\bigl(1-e^{-\imath\Delta T\omega_{j}}\bigr)c\Bigr)\Bigr]^{-1}\\
&=\Bigl[\kappa_{\mathrm{d}}\left(z\right)+\bigl(1-\cos\left(\Delta T\omega_{j}\right)\bigr)c\Bigr]^{-1}\leq
\left[\kappa_{\mathrm{d}}\left(z\right)\right]^{-1}\leq\frac{1}{c_1}.
\end{aligned}
\end{equation}
Due to the Lipschitz condition \eqref{eq:Lipschitz_kappa_d} 
and using \eqref{eq:bdd_inv} in expression \eqref{eq:estim} we directly obtain the estimate \eqref{eq:cond_Newton} with $\delta_0={L_2}/{c_1}$.  
\end{proof}
\end{thm}

%

\IKRrr{Having the constant $\delta_0$ from the Lipschitz condition \eqref{eq:cond_Newton} by the Theorem~\ref{thm:Newton} the convergence of the simplified Newton 
method \eqref{eq:Newton} is then implied,} provided that the bound \eqref{eq:cond_Newton_12} holds 
for {the} selected initial approximation $\mathbf{U}^{(0)}$. 

\IKRrr{From the Theorem~\ref{thm:conv_kappad} we have seen that the convergence result of the 
simplified Newton iteration \eqref{eq:Newton} is derived directly from the properties 
of the nonlinearity of the problem.} However, one should be aware of the fact that 
it {may not always be possible} to find an initial guess $\mathbf{U}^{(0)}$ 
which {is} at the same time {close} enough {to the solution} for convergence of the Newton 
iteration and has the form \eqref{eq:Ub} for applicability of the MH solver 
to system \eqref{eq:Newton_init0} with constant matrix $\mathbf{G}_{\mathrm{d}}$.  

\begin{rem}
\label{rem:diffusion_kappa}
	We note that 
	\IKR{an extension of our convergence results} to the $d$-dimensional problem \eqref{eq:ode_nlin}, 
	obtained from a spatial discretization of a PDE,
	could be based on the operator theory provided in \cite{Pechstein_2004aa}. 
	In this case one would need to consider the properties not only 
	of \IKR{the} function $\kappa$ but also of a nonlinear 
	operator $A\!:V\mapsto V^*,$ 
	which, e.g., for a diffusion equation with nonlinear diffusion coefficient $\kappa(|\nabla u|)$ is defined as
	\begin{equation}
	\big\langle A(u),v \big\rangle_{V^*,V} = \bigl(\kappa(|\nabla u|)\nabla u,\nabla v\bigr)_H 
	\quad\forall u,v\in V,
	\end{equation}
	with a reflexive Banach space $V$ and a Hilbert space $H,$ e.g., $V=H_0^1(\Omega)$ and $H=L^2(\Omega)$ on domain $\Omega.$ 
	\IKRrr{Nevertheless, the necessary conditions for the operator $A$ can be derived directly from the properties of $\kappa$, 
	which are exactly those imposed in the Assumption~\ref{thm:kappa_ass}.}
	We refer to \cite[Corollary 2.13, Theorem 4.2 and Theorem 4.3]{Pechstein_2004aa} for the details.
\end{rem}

\subsubsection{\IKRrr{Numerical study of a model problem}} 
	\begin{figure}[t]
	\centering
	\mbox{
	\includegraphics[width=0.46\linewidth]{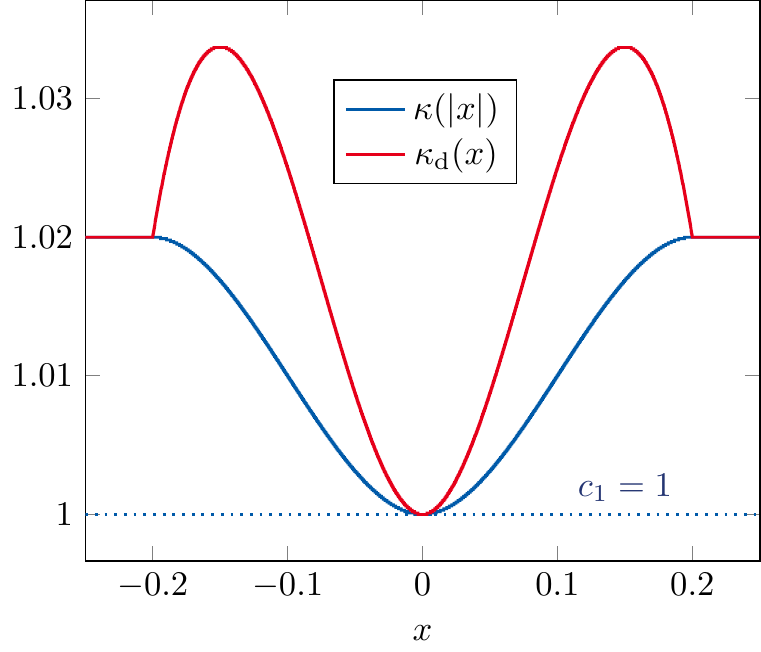}
	\hspace*{0.05\linewidth}
	\includegraphics[width=0.46\linewidth]{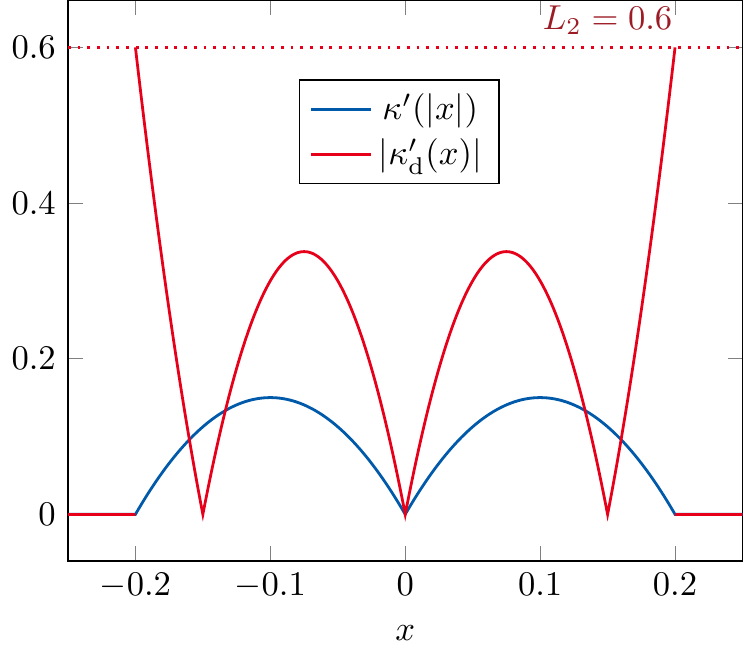}
	}
	\caption[]{\IKRrr{Nonlinearity of the model problem. Left: $\kappa$ and $\kappa_{\mathrm{d}}$. Right: $\kappa^{\prime}$ and $\kappa_{\mathrm{d}}^{\prime}$.}}
	\label{fig:nonlinearity}
	\end{figure}
\IKRrr{
We now discuss convergence of the proposed simplified Newton method applied to a model 1D problem \eqref{eq:ode_1d} 
with $m=10^{-1},$ RHS $j(t)=10^{-3}\sin(2\pi t/T),$ $T=0.02,$ and nonlinearity $\kappa$ defined by a piecewise polynomial 
\begin{equation}
\label{eq:kappa_poly}
\kappa(|x|)=
\begin{cases}
-5|x|^3+1.5|x|^2+1,& 0\leq |x|< 0.1,\\
-5(|x|-0.1)^3+0.15(|x|-0.1)+1.01,& 0.1\leq |x|<0.2,\\
1.02,& 0.2\leq |x|.
\end{cases}
\end{equation}
Such a problem can be interpreted as a nonlinear RL-circuit model with unknown magnetic flux $u$ as in \cite{Gander_2019aa}. 
The nonlinear functions $\kappa$ and $\kappa_{\mathrm{d}}$ are depicted in Figure~\ref{fig:nonlinearity} (left). 
This shows that the lower bound for both $\kappa(|x|)$ and $\kappa_{\mathrm{d}}(x)$, $x\in\mathbb{R}$ is $c_1=1$.  
Figure~\ref{fig:nonlinearity} (right) depicts the corresponding derivatives $\kappa^{\prime}(|x|)$ and $|\kappa_{\mathrm{d}}^{\prime}(x)|,$ $x\in\mathbb{R}.$ 
Based on the observed characteristics of $\kappa$ one can clearly deduce that the Assumption~\ref{thm:kappa_ass} holds. 
Then by Lemma~\ref{thm:kappa_bdd_Lipschitz} the Lipschitz condition 
\eqref{eq:Lipschitz_kappa_d} is satisfied for $\kappa_{\mathrm{d}}$. 
We estimate the Lipschitz constant $L_2$ as the upper bound for $|\kappa_{\mathrm{d}}^{\prime}(x)|\leq L_2=0.6.$ 
Based on Theorem~\eqref{thm:conv_kappad} we have the estimate for $\delta_0=L_2/c_1=6.$
}

\begin{figure}[t]
	\centering
	\mbox{
	\includegraphics[width=0.46\linewidth]{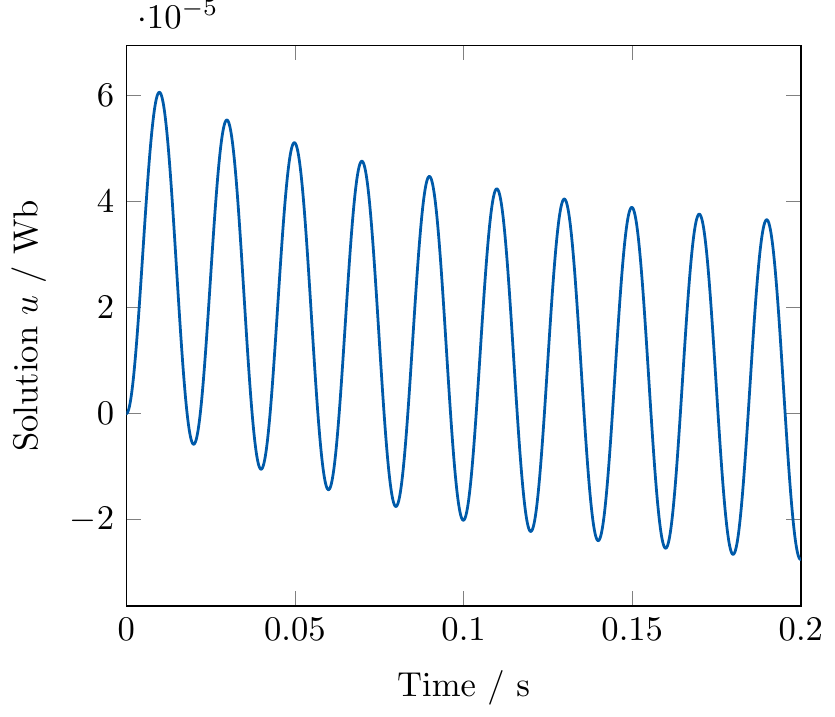}
	\hspace*{0.05\linewidth}
	\includegraphics[width=0.46\linewidth]{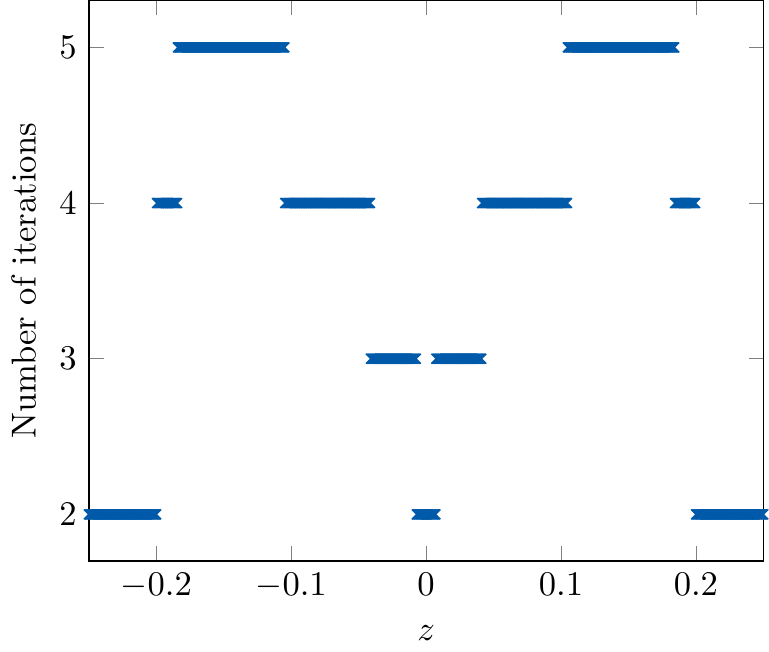}
	}
	\caption[]{\IKRrr{Left: sequential time-stepping over $k^*=10$ periods until the steady state. Right: number of simplified Newton iterations for different choices of the initial approximation $\bu^{(0)}=[z,z,\dots,z]^{\!\top}$.}}
	\label{fig:conv_LR}
	\end{figure}

\IKRrr{The periodic steady-state solution of \eqref{eq:ode_1d} with $\kappa$ from \eqref{eq:kappa_poly} is depicted in Figure~\ref{fig:conv_LR} (left). It is obtained after sequential calculation with step size $\delta T=10^{-5}$ starting from the initial value $u_{\text{init}}=0$ at $t=0$ over $k^*=10$ periods when the tolerance defined in \eqref{eq:error_tol_seq} is reached. We now describe performance of the simplified Newton algorithm (within TP MH setting) for the 1D model problem 
depending on the choice of the initial approximation $\bu^{(0)}=[z,z,\dots,z]^{\!\top}{\in D=(-0.25,0.25)^{N}\subset\mathbb{R}^{N}}$ with $N=10.$ 
In Figure~\ref{fig:conv_LR} on the right we show the number of iterations required until convergence with respect to the norm \eqref{eq:error_tol_IT_TP}. We observed that in the neighborhood of $z=0$ as well as for $|z|\geq0.2$ the algorithm converged in $2$ iterations, while for $0.1<|z|<0.2$ the number of iterations was between $3$ and $5$. 
This shows us that for the considered example the simplified Newton method is convergent for all the considered initial approximations $\bu^{(0)}.$ The calculated periodic solution is $\bu^*\approx10^{-5}\cdot[-3.02,-1.8,0.1,1.9,3.07,3.01,1.8,-0.1,-1.96,-3.07]^{\!\top}$ in all cases. }

\begin{figure}[t]
	\centering
	\mbox{
	\includegraphics[width=0.46\linewidth]{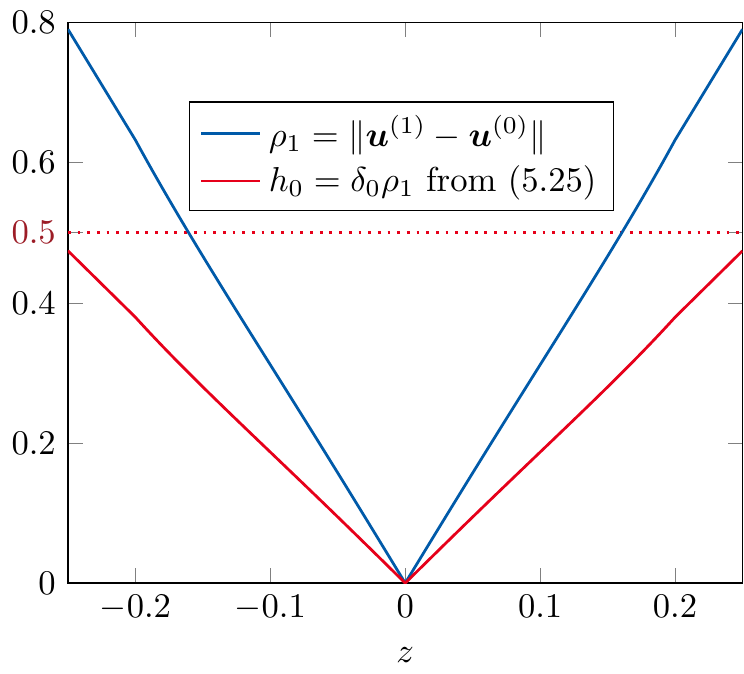}
	}
	\caption[]{\IKRrr{Calculated deviation of $\bu^{(1)}$ form the initial approximation $\bu^{(0)}$ and the constant $h_0$ from \eqref{eq:cond_Newton_12}. 
	}}
	\label{fig:conv_LR_details}
	\end{figure}

\IKRrr{For each $\bu^{(0)}$ we have calculated the values of $h_0=\delta_0\rho_1=0.6\rho_1$ from \eqref{eq:cond_Newton_12} with $\rho_1=\|\bu^{(1)}-\bu^{(0)}\|$ which we depicted in Figure~\ref{fig:conv_LR_details}. We see that condition \eqref{eq:cond_Newton_12} holds, since for all cases $h_0<0.5.$ Then by the Theorem~\ref{thm:Newton} the iterates $\bu^{(s)}$ remain in the ball $\bar{S}\left(\bu^{(0)},\rho\right)$ of radius $\rho$ defined in \eqref{eq:rho_Newton}, which we have also observed in practice.} 



	\section{Numerical experiments for a two-dimensional coaxial cable model}\label{sec:numerics}
	In this section we present the results obtained from application 
	of the proposed algorithm with MH coarse correction 
	to solve the time-periodic problem for 
	a coaxial cable model. Figure~\ref{fig:tube} (left) illustrates a sketch of a two-dimensional (2D) domain 
	\begin{equation}
	\bar{\Omega}=\bar{\Omega}_{\mathrm{Fe}}\cup\bar{\Omega}_{\mathrm{Cu}}\cup\bar{\Omega}_{\mathrm{Air}}.
	\end{equation}
	It represents the cross-section of the cable,
	which consists of steel tube $\Omega_{\mathrm{Fe}},$ 
	conducting wire $\Omega_{\mathrm{Cu}},$ and air gap 
	$\Omega_{\mathrm{Air}}.$
	
	{The} electromagnetic phenomena in the coaxial cable when 
	the inner wire is supplied by a sinusoidal current source {are}  
	mathematically described by the eddy current problem \cite{Jackson_1998aa}. 
	In the following subsection we discuss the mathematical modeling of the 
	underlying physical characteristics and describe the properties of the 
	nonlinearity. 
	
	\subsection{Eddy current problem} 
  A time-periodic eddy current formulation in terms of {the} magnetic vector potential $\vec{A}$ 
	\cite{Jackson_1998aa} searches for unknown $\vec{A}(\vec{x},t),$ with 
	$(\vec{x},t)\in\Omega\times[0,T],$ which solves 
	\begin{align}
	\sigma\left(\vec{x}\right)\partial_t\vec{A}(\vec{x},t)+
	\nabla\times\Bigl(\nu\bigl(\vec{x},|\nabla\times\vec{A}|\bigr)\nabla\times\vec{A}(\vec{x},t)\Bigr) &= \vec{j}(\vec{x},t) &&\text{{in}}\ \Omega\times(0,T),\label{eq:mqs1}    \\
	\vec{A}(\vec{x},0)&=\vec{A}(\vec{x},T), &&\vec{x}\in\Omega,\label{eq:mqs1_init}\\
	{\vec{n}}\times\vec{A}&= 0 &&\text{{on}}\ \partial\Omega{\times[0,T]},\label{eq:mqs1_bdry}
	\end{align}
	where $\partial\Omega$ denotes the boundary of $\Omega$ and $\vec{n}$ is the outward normal vector to $\partial\Omega$.  
	\IKR{The} function $\sigma(\vec{x})\geq 0,$ $\vec{x}\in\Omega$ {describes} \IKR{the} electric conductivity of 
	the materials. It is positive only for $\vec{x}\in\Omega_{\mathrm{Fe}}$ and 
	is equal to zero in $\Omega\setminus{\Omega}_{\mathrm{Fe}}.$ This gives \IKR{a}  
	parabolic-elliptic character to \eqref{eq:mqs1}. The sinusoidal current excitation  
	defines the RHS
	\begin{equation}
	\vec{j}(\vec{x},t)=\mathbbm{1}_{\Omega_{\mathrm{Cu}}}(\vec{x})100\sin(2\pi t/T),\ \ (\vec{x},t)\in\Omega\times(0,T),
	\end{equation}
	where $\mathbbm{1}_{\Omega_{\mathrm{Cu}}}$ denotes the indicator function of subdomain $\Omega_{\mathrm{Cu}}.$ 	
	The input $\vec{j}$ \IKR{is} described using winding functions \cite{Schops_2013aa}.
	{The} magnetic reluctivity $\nu\bigl(\vec{x},|\nabla\times\vec{A}|\bigr)>0$ is defined by the ferromagnetic 
	material in $\Omega_{\mathrm{Fe}}$ and {equals} the reluctivity of vacuum $\nu_0=\SI{10^7/(4\pi)}{\metre\per\henry}$
	in $\Omega\setminus{\Omega}_{\mathrm{Fe}}.$ 
	
	\begin{figure}[t]
	\centering
			\mbox{\includegraphics[width=0.45\linewidth]{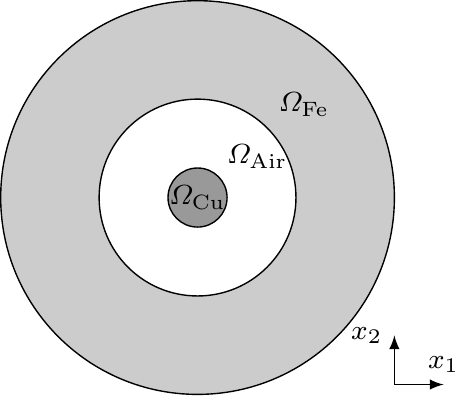}
			\hspace*{0.02\linewidth}
			\includegraphics[width=0.5\linewidth]{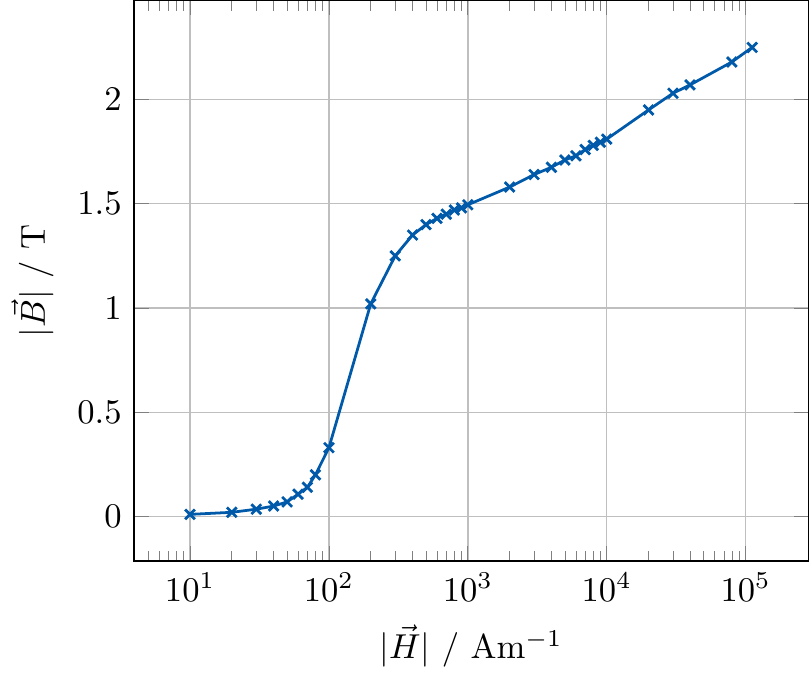}}
	\caption[]{Left: cross-section of a coaxial cable\footnotemark[3]. 
	Computational domain consists of tube $\Omega_{\mathrm{Fe}}$, conducting wire 
	$\Omega_{\mathrm{Cu}}$, and air gap $\Omega_{\mathrm{Air}}$ in-between. Right: 
	magnetization curve of the ferromagnetic material in $\Omega_{\mathrm{Fe}}$ 
	given by measured data points.}
	\label{fig:tube}
	\end{figure}	
	\footnotetext[3]{\url{http://www.femm.info/wiki/tubeexample}}

	{The} reluctivity function in $\Omega_{\mathrm{Fe}}$ is determined  
	by {the} magnetization curve $b:\mathbb{R}^+_0\to\mathbb{R}^+_0$, shown in {the} Figure~\ref{fig:tube} (right). 
	It is commonly obtained from measured tuples $(|\vec{B}|,|\vec{H}|)$ so that $|\vec{B}|=b\bigl(|\vec{H}|\bigr),$ 
	where $|\vec{B}|$ denotes the magnitude of the magnetic flux density $\vec{B}$ and $|\vec{H}|$ is the   magnitude of the 
	magnetic field intensity $\vec{H}$. To obtain a continuously differentiable curve  
	interpolation of the measurement data is performed, e.g., using classical splines or 
	specific approximation techniques as in \cite{Pechstein_2006aa} \todo{such that the monotonicity is preserved}. 	
	We can then define {the} reluctivity 
	\begin{equation}
	\label{eq:nu_all}
	\nu\bigl(\vec{x},|\vec{B}|\bigr)=
	\begin{cases}
	\nu_0, & \vec{x}\in\Omega\setminus{\Omega}_{\mathrm{Fe}},\\
	\nu(|\vec{B}|), &\vec{x}\in{\Omega}_{\mathrm{Fe}}.
	\end{cases}
	\end{equation}
	For {each} $\vec{x}\in{\Omega}_{\mathrm{Fe}}$ we {can construct} a continuously differentiable 
	reluctivity function $\nu:\mathbb{R}^+_0\to\mathbb{R}^+_0$ 
	from the magnetization curve as \cite{Heise_1994aa}
	\begin{equation}
	\label{eq:nu}
		\nu(s)=\begin{cases}
		\displaystyle\frac{b^{-1}(s)}{s},& s\in\mathbb{R}^+,\\
		(b^{-1})'(0), & s=0.
		\end{cases}
	\end{equation}
	Based on the physical nature of magnetization curve $b$ we can naturally impose the following 
	assumptions on \IKR{the} reluctivity $\nu$ \cite{Heise_1994aa}.
	\begin{ass}[Conditions on the magnetic reluctivity]
	\label{ass:nu}
	Let for $\nu:\mathbb{R}^+_0\to\mathbb{R}^+_0$ it holds:
	\begin{itemize}
	\item[$\bullet$] $\nu(t)\geq c_1>0,$ $t\in\mathbb{R}_0^+,$
	\item[$\bullet$] function $\IKRrr{f}:t\rightarrow\nu(t)t$ is strongly monotone with monotonicity constant $c_1,$ see \eqref{eq:str_monotone},
	\item[$\bullet$] $\nu\in C^1(\mathbb{R}_0^+)$ and $\lim\limits_{t\to\infty}\nu^{\prime}(t)=0,$ 
	\item[$\bullet$] function $\IKRrr{g}:t\rightarrow\nu^{\prime}(t)t$ satisfies Lipschitz condition \eqref{eq:Lipschitz_cont}.
\end{itemize}
\end{ass}
\IKRrr{We note that the conditions on $\nu$ are the same as those imposed on the function $\kappa$ in Assumption~\ref{thm:kappa_ass}.  Analagously to the Remark~\ref{rem:diffusion_kappa} 
these properties of the nonlinear reluctivity function $\nu$ define 
conditions on the nonlinear operator in \eqref{eq:mqs1}, necessary  
for convergence of the simplified Newton algorithm \eqref{eq:Newton_init0}.} 

For \IKR{the} numerical solution of \eqref{eq:mqs1}-\eqref{eq:mqs1_bdry} the spatial discretization is performed, e.g., 
using FEM with linear shape functions \cite{Brenner_2008aa}. This leads to a time-dependent system of 
differential-algebraic equations (DAEs) \cite{Lamour_2013aa} \IKRrr{(due to the presence of the non-conducting domain $\Omega\setminus{\Omega}_{\mathrm{Fe}}$)}, which together with the periodicity constraint 
has \IKR{the} form of \eqref{eq:ode_nlin}, 
with \IKR{a} (singular) mass matrix $\mathbf{M}\in\mathbb{R}^{d\times d},$ nonlinear stiffness matrix 
$\mathbf{K}(\cdot):\mathbb{R}^d\to\mathbb{R}^{d\times d}$ and unknown $\bu:[0,T]\to\mathbb{R}^d.$ 

We consider $T=0.02$ s to be the period in the eddy current formulation \eqref{eq:mqs1}-\eqref{eq:mqs1_bdry}. 
In our numerical simulations a FEM discretization of \IKR{the} domain $\Omega$ with $d=2269$ degrees of freedom is exploited. 
The time integration is performed with the implicit Euler method. \todo{The DAE is of index-1, therefore in this case it 
does not need any extra care, e.g., with respect to consistent initial values \cite{Schops_2018aa}.} \IKRr{The nonlinear algebraic system 
at each time step is solved using the Newton method.} The time step size $\delta T=10^{-5}$ s is used for the sequential calculation, for TP discretization \eqref{eq:sys_explicit_nlin}, 
and also as a step size on the fine level of the considered parallel-in-time methods. The coarse propagator solves only one time 
step per window $[T_{n-1},T_n],$ $n=1,\dots,N,$ thereby using step size $\Delta T = T/N.$ 

We \IKR{now} illustrate {the} performance of the MH approach 
for linear and nonlinear models of the coaxial cable and compare it with existing standard and 
parallel-in-time approaches. \IKRrr{We discuss the performance firstly in terms of effective linear system solved. This number corresponds to the solver calls executed sequentially which, in case of the parallel computations, is the maximum number of linear system solutions performed by a single core. This is a convenient and easily reproducible measure because of its independence with respect to the software implementation. However, it disregards aspects like communication costs and thus also wallclock times are compared. }

\subsection{Linear model} \label{subsec:numerics_linear} 
	\begin{figure}[t]
	\centering 
	\mbox{	
		\includegraphics[width=0.475\linewidth]{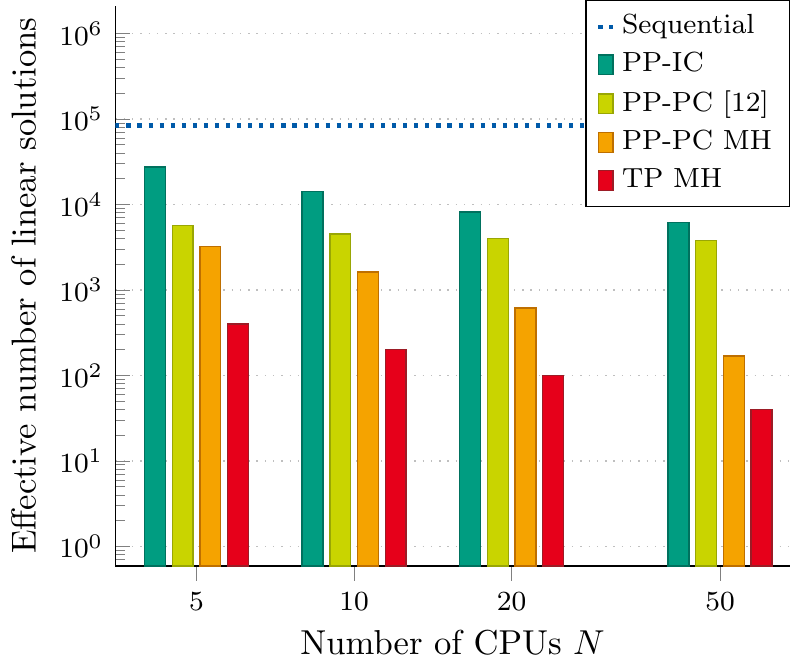}
		\hspace*{0.02\linewidth}
\includegraphics[width=0.475\linewidth]{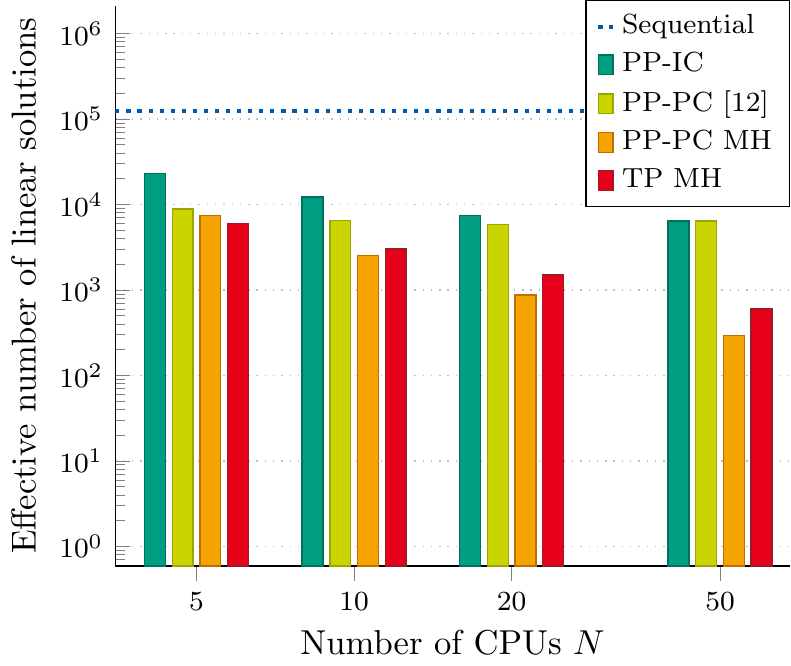}
}
	\caption{Comparison of the computational costs for different approaches, applied to the coaxial cable model: strong scaling. \IKRrr{Left: linear model. Right: nonlinear model.} \IKR{Comparing bars of different colors for fixed $N$ shows the \IKRrr{performance gain due to the methods}, comparing bars of the same color corresponds to the \IKRrr{gain due to parallelism}.}}
	\label{fig:res_lin}
	\end{figure}
	Let us assume that the reluctivity function $\nu$ in \eqref{eq:mqs1} is linear, 
	i.e., it does not depend on {the} solution $\vec{A}.$ We therefore deal with problem 
	\eqref{eq:ode_lin}, where {the} matrix $\mathbf{K}$ is constant. 
	One could then construct a linear (PP-PC) 
	system of form \eqref{eq:PP_PC_sys_AVG_lin} and directly transform it into the frequency 
	domain \eqref{eq:PPPC_freq_domain} \cite{Kulchytska-Ruchka_2019ac}, thereby solving $N$ 
	separate $d$-dimensional systems in parallel. 
	
	\IKR{The} computational costs of several solution approaches are compared \IKRrr{for the linear model} in terms of the effective 
	number of required linear system solutions. 
	%
	\IKR{Figure~\ref{fig:res_lin} on the left} illustrates how the costs of the parallel algorithms \IKRrr{applied to the linear model} scale 
	when the number of cores $N$ increases, while the precision of 
	the fine solution remains unchanged (strong scaling). 
	The classical time-stepping approach, 
	when {the} solution of the initial value problem is performed sequentially 
	starting from zero initial condition, required time integration over $\IKR{k^*=42}$ 
	periods \IKRrr{(which corresponds to $84\,000$ time steps and solution of linear systems)} until the \IKR{error \eqref{eq:error_tol_seq} of the obtained} periodic 
	steady-state solution reached the tolerance\IKR{s} specified in \eqref{eq:error_tol}. \IKRrr{We consider this sequential time-stepping solution as a benchmark to compare performance of the other applied approaches.}
	
	The application of PP-IC \eqref{eq:PP-IC1}-\eqref{eq:PP-IC2} \IKRrr{reduced the number of effective linear system solutions} 
	by \IKR{a} factor \IKR{of} $14$ \IKR{compared to} the sequential calculation when using $N=50$ cores and by a factor of $3$ with $N=5$ CPUs. 
	{The} PP-PC method \IKRrr{with the Jacobi-like iteration scheme} \eqref{eq:PP-PC_sys_fxd_pt} calculates the 
	periodic solution \IKRrr{effectively solving $15$ times less linear systems} than the classical time stepping already when exploiting only $N=5$ processing units. 
	However, the disadvantage of this approach is that the computational costs stagnate and there is no further acceleration when 
	applying additional computational power. 
	
	In contrast to this, the MH approach applied to the PP-PC system \eqref{eq:PP_PC_sys_AVG_lin} as in \cite{Kulchytska-Ruchka_2019ac}  
	delivers the steady-state solution \IKRrr{effectively solving $500$ times less linear systems when using $N=50$ and $27$ times less linear systems when using $N=5$ cores,} compared to the standard time integration. Finally, {the direct} MH solution of the TP system 
	\eqref{eq:sys_explicit_nlin} {is} performed on the fine grid (with step size $\delta T=10^{-5}$ s on $[0,0.02]$ s). {It} requires $N_{\mathrm{f}}=2\,000$ 
	separate linear system solves. Clearly, in the linear case the \IKRrr{performance gain} obtained by this approach scales \IKRrr{perfectly} linear as $N$ grows. 
	In particular, when distributing the workload among $N=5$ CPUs the \IKRrr{computational efforts of the steady-state calculations in terms of the effective linear system solves are reduced} by a factor of $210$ and for $N=50$ by a factor of $2\,100.$
	Therefore, as expected the TP solution with the MH approach \cite{Biro_2006aa} {yields optimal scaling} for 
	the solution of the linear time-periodic eddy current problem.

\subsection{Nonlinear model}\label{subsec:numerics_nonlinear} 
	We now present {the} performance of the iterative method introduced in Section~\ref{section:Parareal_new} 
	for parallel-in-time solution of {a} nonlinear time-periodic {problem} \eqref{eq:mqs1}-\eqref{eq:mqs1_bdry}. 
	In this case the nonlinear time-periodic system \eqref{eq:PP_PC_sys_explicit_nlin} is solved with the simplified Newton 
	method \eqref{eq:Newton}, accompanied by the frequency domain solution at each {Newton} iteration. The initial {guess in 
	\eqref{eq:U0} for the simplified Newton iteration uses $\mathbf{Z}=\mathbf{0}$, i.e.,
	\begin{equation}
	\mathbf{U}^{(0)}=\IKR{\mathbf{U}^{(k+1,0)}}:=\left[\bigl(\mathbf{b}_N^{(k)}\bigr)^{\!\top},\bigl(\mathbf{b}_1^{(k)}\bigr)^{\!\top}\dots,\bigl(\mathbf{b}_{N-1}^{(k)}\bigr)^{\!\top}\right]^{\!\top}
	\end{equation} 
	is chosen at PP-PC iteration $k+1$.} \IKRrr{Almost identical results have been obtained for the choice
	$$ \mathbf{Z}^{(k)} = \frac{1}{N}\sum\limits_{n=0}^{N-1}\mathbf{U}_n^{(k)}, $$
    i.e., the average among the $N$ time steps at PP-PC iteration $k$. Since such choices are typically problem-specific, there is no further discussion of `sufficiently good' initial guesses here. The interested reader is referred to, e.g., \cite{Deuflhard_2004aa}}
	
\begin{figure}[t]
\centering
\mbox{	
\includegraphics[width=0.475\linewidth]{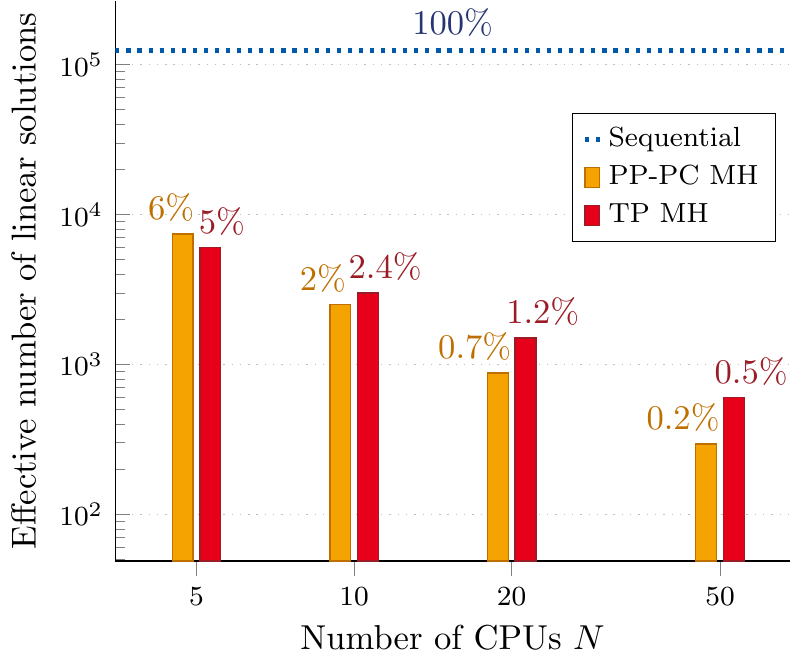}
\hspace*{0.02\linewidth}
\includegraphics[width=0.475\linewidth]{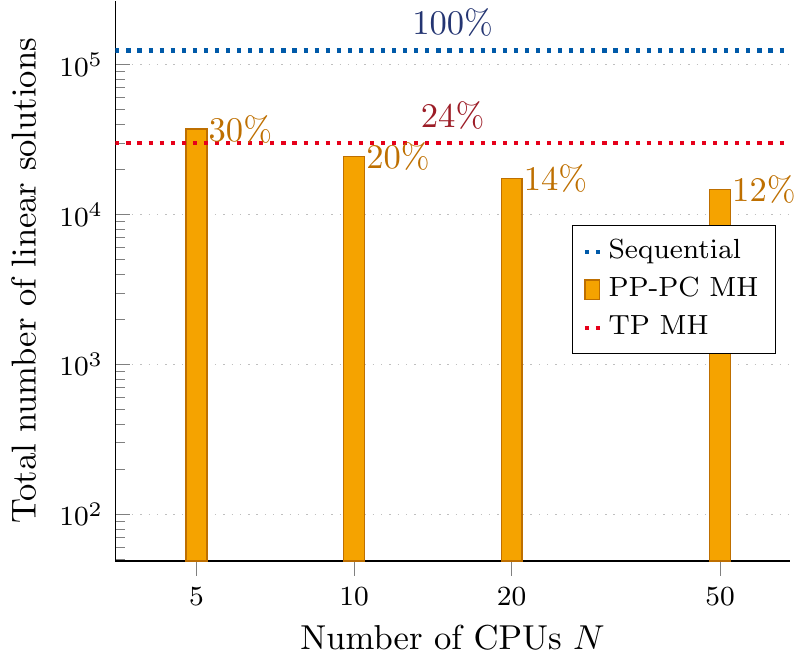}
}
	\caption{\IKRrr{Number of linear systems solved within PP-PC MH and TP MH for the nonlinear coaxial cable model and its comparison to the sequential time stepping. Left: effective number of linear solutions (maximum among CPUs). Right: total number of linear solutions (sum for all CPUs).}}
	\label{fig:res_nlin_MH}
\end{figure}
	
	As in Section~\ref{subsec:numerics_linear}, Figure~\ref{fig:res_lin} (on the right) illustrates the computational efforts of all the previously described 
	methods to obtain the periodic steady-state solution for the nonlinear coaxial cable model. 
	Sequential time stepping starting from zero initial value reached the steady-state solution, periodic  
	up to the tolerance\IKR{s} {given} in \eqref{eq:error_tol} \IKR{using the norm \eqref{eq:error_tol_seq}}, after calculation over $\IKR{k^*=31}$ periods, \IKRrr{which corresponded to solution of $124\,000$ linear systems.} PP-IC computed the periodic solution 
	\IKR{\IKRrr{effectively} solving $19$ times less linear systems when using $N=50$ cores and $5$ times less systems when $N=5.$ The reduction of the effective number of solved linear systems} using PP-PC \eqref{eq:PP-PC_sys_fxd_pt} amounts to a factor of $14$ with $N=5$ CPUs and to a factor of $19$ with $N=50,$ as within PP-IC. As previously seen, the costs of the PP-PC solution 
	do not significantly decrease with the growth of $N$ and even increase slightly, when comparing the results for $N=20$ and $N=50.$ 
	
	On the other hand, the simplified Newton iteration \eqref{eq:Newton} for PP-PC \eqref{eq:PP_PC_sys_explicit_nlin} with MH correction 
	on the coarse grid \IKRrr{effectively solves} about $16$ \IKRrr{times less linear systems} for $N=5$ and $420$ \IKRrr{times less} for $N=50.$ Now, \IKRrr{in contrast} to the linear 
	case, the TP MH solution (with $\bu^{(0)}=\mathbf{0}$ in \eqref{eq:U0_TP}) performs worse than PP-PC MH for almost all the considered values of $N$ except $N=5.$ More specifically, for $10$, 
	$20$, and $50$ cores the \IKRrr{number of effective linear solutions within} TP MH \IKRrr{is bigger} than that of PP-PC MH, e.g., reduction factor {is} equal 
	to $206$ \IKR{(versus factor $420$ of PP-PC MH) when using $N=50.$ In case of $N=5,$ $20$ times less systems were solved effectively with TP MH  which is \IKRrr{a better result} than that of 
	PP-PC MH (where factor $16$ is observed).}	
	\IKRrr{
	In Figure~\ref{fig:res_nlin_MH} on the left we see that the percentages of the corresponding calculations with respect to the sequential solution are considerably small (less than $1\%$ for $N=50$ and only $5$-$6\%$ for $N=5$). 
	}
	
	\IKRrr{The Figure~\ref{fig:res_nlin_MH} (right) illustrates the total number of linear systems solved within PP-PC MH and TP MH, i.e., the numbers of linear systems solved by each CPU are added. We see that TP MH always solves the same number of linear systems, since increasing the number of cores does not influence the convergence. The method converged in $15$ iterations, each calculating $2\,000$ steps, which in total gives $15\times2\,000=30\,000$ systems solves. In contrast to this, the reduction of the computational costs within PP-PC MH comes not only from distributing the workload among CPUs but also from its faster convergence with bigger $N$ due to the increasing accuracy of the coarse solver. Therefore, even when the code is not truly parallelized, both approaches solve much less linear systems compared to the classical sequential solution ($4$ times less for TP MH and from $3$ up to $8$ times less systems for PP-PC MH).}


\begin{figure}[t]
\centering
\mbox{	
\includegraphics[width=0.485\linewidth]{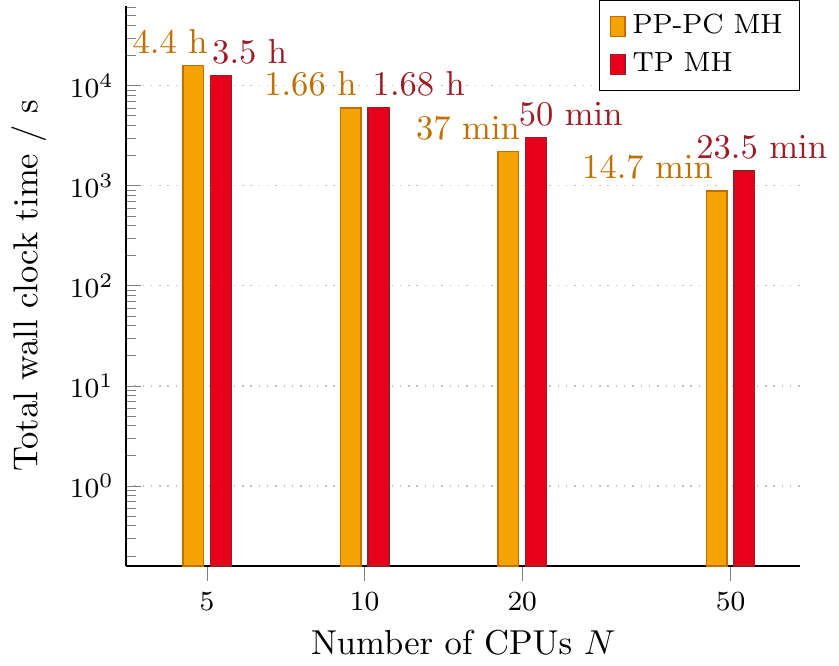} 
\hspace*{0.02\linewidth}
\vspace*{1em}
\includegraphics[width=0.465\linewidth]{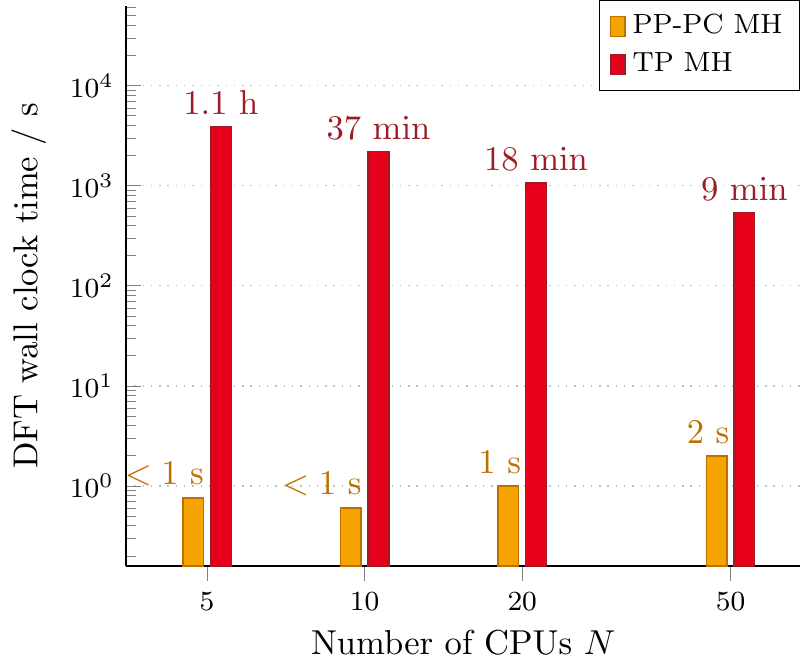}
}
\caption{\IKRrr{Measurement of the wall clock time of the simplified Newton MH approach, applied to the PP-PC and the TP systems for the nonlinear coaxial cable model. Left: total computational time. Right: time spent on calculation of the DFT.} }
	\label{fig:res_nlin_timings}
\end{figure}
\IKRrr{Finally, we compare the actual computational time of the simplified Newton approach, applied to the PP-PC and the TP systems, followed by the MH correction. 
As before, the wall clock time measurements illustrate better performance of PP-PC MH in case of $N=10$, $20$ and $50$, and of TP MH when $N=5$ (see Figure~\ref{fig:res_nlin_timings} on the left). For $N=5$ PP-PC MH took $4.4$ hours, while TP MH required $3.5$ hours. In case of $N=50$ calculation with PP-PC MH lasted $14.7$ minutes, whereas TP MH took $23.5$ minutes. In comparison, the sequential time stepping over one period was executed in $2.15$ hours, which would extend to a simulation of about $66.7$ hours (almost $3$ days) calculating over $31$ periods until the steady state is reached. Therefore, for PP-PC MH the speed up factor is $15$ for $N=5$, $40$ for $N=10$, $108$ for $N=20$, and $271.5$ for $N=50.$ In its turn, TP MH accelerates the sequential time stepping by a factor of: $19$ for $N=5$, $40$ for $N=10$, $79$ for $N=20$, and $170$ for $N=50.$ The corresponding percentages with respect to the sequential execution time are presented in Table~\ref{tab:Comparison}.}

\IKRrr{
On the right in Figure~\ref{fig:res_nlin_timings} the time needed to calculate the Fourier transform together with its inverse is illustrated. We see that within TP MH the DFT occupies approximately a third part of the total computational time for all $N$ (from $33\%$ to $37\%$), while the DFT cost of PP-PC MH is extremely small (less than $1\%$ of the total execution time.) Besides, we see that for TP MH the time for DFT decreases with the increase of $N$, as expected due to the better parallelization capability. The opposite happens within PP-PC MH: although the cost of DFT still remains tiny it grows together with $N$, since the size of the DFT matrix $\mathbf{F}$ \eqref{eq:Fpq} increases. Finally, we would like to comment that the communicatioin costs of the both approaches are within the range of $14$ s and $48$ s for the considered values of $N.$ }

\begin{figure}[t]
	\centering
\mbox{	
\includegraphics[width=0.475\linewidth]{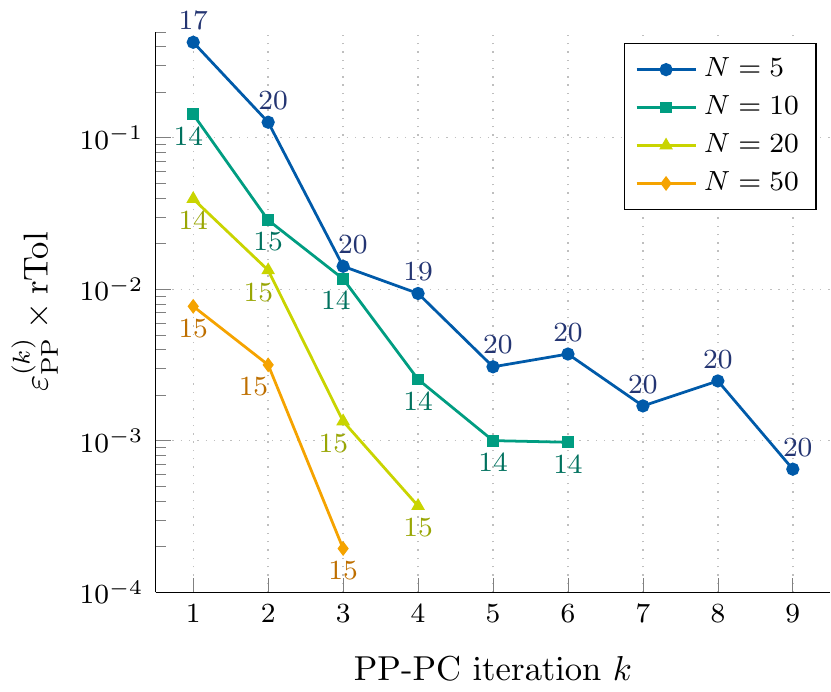}
\hspace*{0.02\linewidth}
\vspace*{1em}
\includegraphics[width=0.475\linewidth]{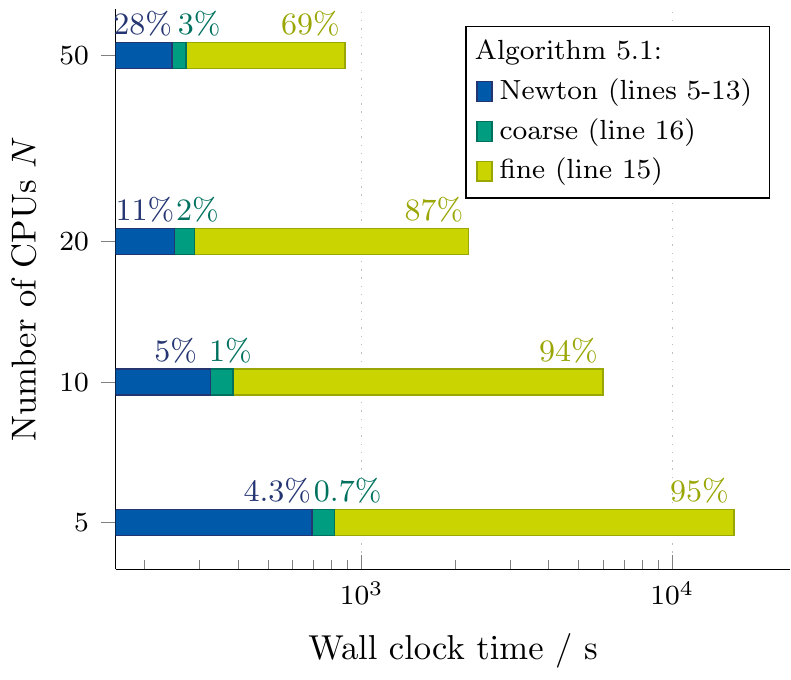}
}
\caption{\IKRrr{Performance of PP-PC MH for the nonlinear coaxial cable model. Left: error and number of the simplified Newton iterations at each PP-PC iteration $k$. Right: contribution of different components of the algorithm to the total computational time.}}
	\label{fig:res_nlin_pppc_mh}
\end{figure}

\begin{table}
		\caption{\IKRrr{Percentage of the wall clock time of PP-PC MH and TP MH with respect to the sequential time stepping.} }
		\centering
		\begin{tabular}{@{}lcccc@{}}
			\toprule
			& \IKRrr{$N=5$} & \IKRrr{$N=10$} & \IKRrr{$N=20$} & \IKRrr{$N=50$} \\
			\midrule
			\IKRrr{Sequential} & \multicolumn{4}{c}{\IKRrr{$100\%$}}\\
			\IKRrr{PP-PC MH}  & \IKRrr{$6.6\%$} & \IKRrr{$2.5\%$} & \IKRrr{$0.9\%$} & \IKRrr{$0.4\%$}\\
			\IKRrr{TP MH} & \IKRrr{$5.3\%$} & \IKRrr{$2.5\%$} & \IKRrr{$1.3\%$} & \IKRrr{$0.6\%$}\\
			\midrule
			\bottomrule
		\end{tabular}		
		\label{tab:Comparison}
	\end{table}

\IKRrr{On the left in Figure~\ref{fig:res_nlin_pppc_mh} we show the convergence of the PP-PC MH approach for the strong scaling. As expected for a Parareal-based method, the more subintervals $N$ are used, the less iterations are required until convergence of PP-PC due to the increasing precision of the coarse solver. This is also represented on the right in Figure~\ref{fig:res_nlin_pppc_mh}, where we illustrate the percentage of the principal components of PP-PC MH regarding the total calculation time. In particular, each iteration of PP-PC includes the parallel fine and coarse solutions and the simplified Newton iterations.  As a result of the decreasing number of PP-PC iterations with a bigger $N$ the timings of each part of the algorithm decreases. The most time-consuming part is the parallelized fine solution, whose contribution amounts to $95\%$ of the total time when $N=5$ and to $69\%$ for $N=50$. The calculation timings of the coarse solution and of simplified Newton also decrease, while their parts of the total time are respectively $4.3\%$ and $0.7\%$ for $N=5$; and $28\%$ and $3\%$ when $N=50.$}

		\IKR{We note that in both linear and nonlinear cases the errors \eqref{eq:error_tol_PP}, 
		\eqref{eq:error_tol_IT}-\eqref{eq:error_tol_seq} were respectively calculated for the 
		corresponding solution approaches. We have additionally performed comparison of the obtained 
		results to the benchmark sequential solution using the norm \eqref{eq:error_tol} at each 
		fine time step and calculating the maximum among the values. For all the considered methods 
		the calculated solutions were close to the benchmark solution in the norm \eqref{eq:error_tol} up to the tolerances 
		of $\text{rTol}=2.5\cdot 10^{-2}$ and $\text{aTol}=2.5\cdot 10^{-5}.$}	

	\section{\IKRr{Application to a three-dimensional transformer model}}\label{sec:numerics_transfo}
\begin{figure}[t]
	\centering
			\mbox{
			\includegraphics[width=0.431\linewidth]{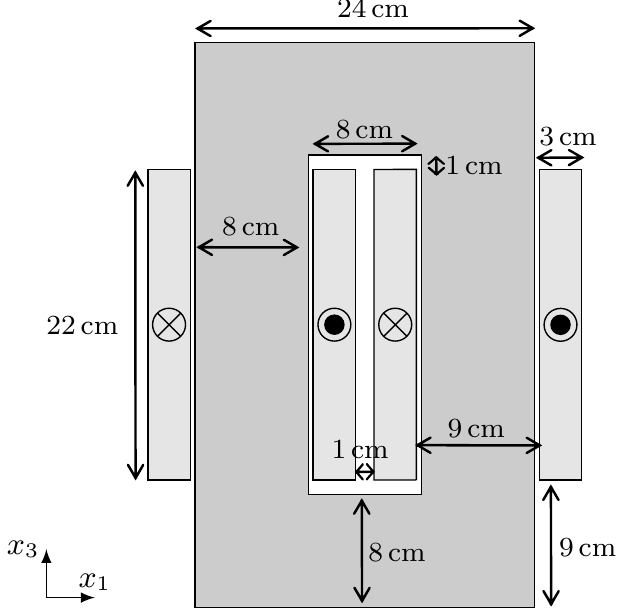}
			\hspace*{0.075\linewidth}
			\includegraphics[width=0.3\linewidth]{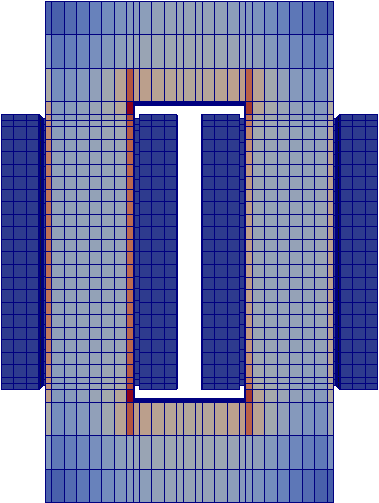}
			\hspace*{0.005\linewidth}
			\includegraphics[width=0.07\linewidth]{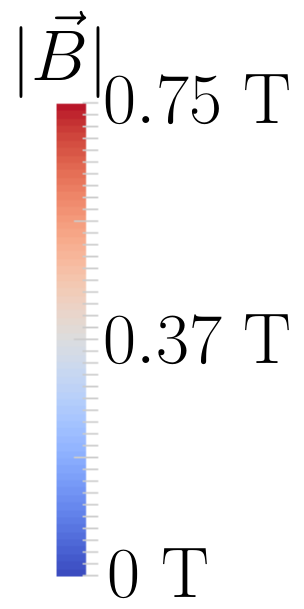}
			}
	\caption[]{Left: cross-section of the 3D FIT model of the transformer. 
			   Right: magnetic flux density distribution at $t=0.02$ s.}
	\label{fig:transfo}
	\end{figure}

	We now apply the proposed algorithm with the MH correction 
	to a three-dimensional (3D) model of a transformer, whose cross-section and dimensions are illustrated in 
	Figure~\ref{fig:transfo} on the left. The computational domain $\Omega$ consists of a  
	steel core and two copper coils, surrounded by air with prescribed homogeneous Dirichlet conditions on the outer boundaries. 
	The total depth of the coils in the dimension $x_2$ is $10$ cm, while the width of the core is $3$ cm. 
	The two coils have $358$ and $206$ windings wound around the core, respectively. 
	
	The nonlinear magnetic reluctivity $\nu$ is given in the ferromagnetic material of the core 
	by the Brauer's curve \cite{Brauer_1981aa}
	\begin{equation}
	\nu\bigl(\vec{x},|\vec{B}|\bigr)=k_1\exp(k_2|\vec{B}|^2)+k_3,
	\end{equation}
	with parameters $k_1=0.3774,$ $k_2=2.970,$ and $k_3=388.33.$ 
	The reluctivities in the air and the copper regions are given by the reluctivity of vacuum $\nu_0.$	
	
	The electric conductivity $\sigma(\vec{x})$ is nonzero only in the steel part and is equal to $\SI{5\cdot 10^{5}}{\siemens\per\metre}$. 
	The coils provide the sinusoidal current excitation 
	\begin{equation}
	\vec{j}(\vec{x},t)=\vec{\chi}_1(\vec{x})10\sin(2\pi t/T)+\vec{\chi}_2(\vec{x})2\sin(2\pi t/T),\ \ (\vec{x},t)\in\Omega\times(0,T),
	\end{equation}
	where $\vec{\chi}_i(\vec{x}),$ $i=1,2$ denote the winding functions of the two coils \cite{Schops_2013aa} and $T=0.02$ s.  			

	Using these input data we consider the eddy current problem 
	of the form \eqref{eq:mqs1}-\eqref{eq:mqs1_bdry} here. 
The spatial discretization is performed using the finite integration technique (FIT) \cite{Weiland_1996aa} with $d=48\,417$ degrees of freedom, which leads to a time-periodic problem of the form \eqref{eq:ode_nlin}. 
\IKRrr{Although the system size is considerably bigger than that of a coaxial cable model from Section~\ref{sec:numerics}, we note that} the considered grid illustrated in Figure~\ref{fig:transfo} on the right is visually relatively coarse.

\begin{figure}[t]
	\centering
	\mbox{
	\includegraphics[width=0.45\linewidth]{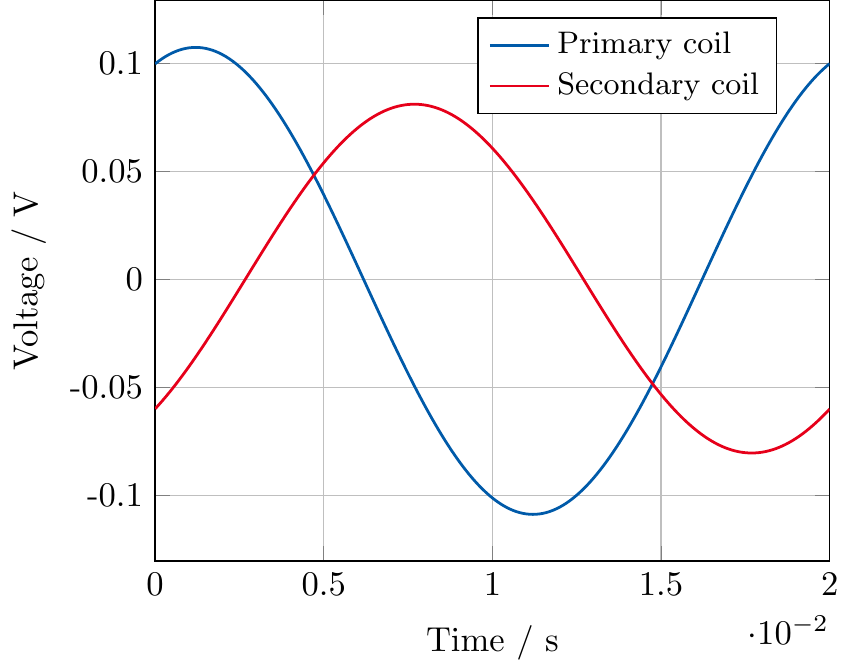}
	\hspace*{0.05\linewidth}
	\includegraphics[width=0.35\linewidth]{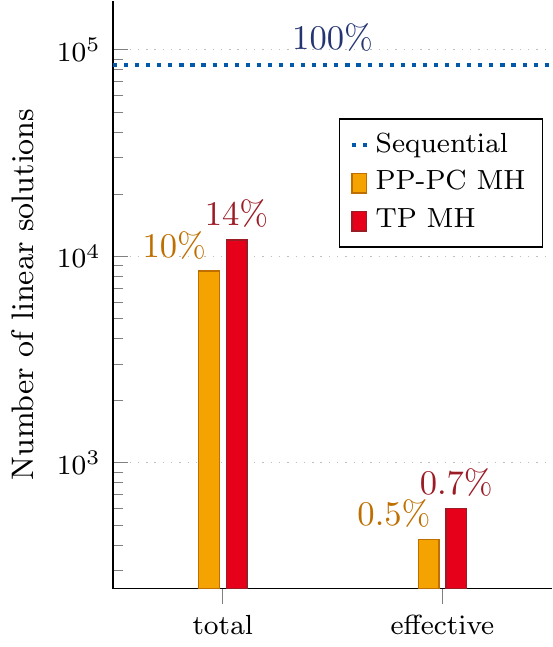}
	}
	\caption[]{Left: induced time-periodic voltages in the primary and secondary coils over the period $[0,0.02]$ s. \IKRrr{Right: total and effective number of linear systems solved within PP-PC MH and TP MH for $N=20$ in contrast to sequential time stepping.}}
	\label{fig:transfo_vltg}
	\end{figure}

The time integration is achieved with the implicit Euler method. The nonlinear system of equations at each time step 
is solved using the successive substitution method, i.e., for fixed $n\geq1$ and for $s=0,1,\dots$ the solution $\bu_n^{(s+1)}$ of the linear system
\begin{equation}
 \left[\frac{1}{\delta T}\mathbf{M}+ \mathbf{K}\bigl(\bu_n^{(s)}\bigr)\right]\bu_n^{(s+1)} =\frac{1}{\delta T}\mathbf{M}\bu_{n-1}+\bj(t_n)
 \label{discr_eddy_current_TP_succ_subst}
\end{equation}
is calculated. The fine step size $\delta T=10^{-5}$ s and the coarse step size $\Delta T = 10^{-3}$ s are used 
within the (parallel-in-time) time-domain simulations.  
Parallelization among $N=20$ CPUs is exploited within the considered parallelization approaches.
The solutions are calculated up to the tolerances specified in \eqref{eq:error_tol}. 
The resulting magnetic flux density distribution at $t=0.02$ s is depicted in the 
Figure~\ref{fig:transfo} on the right. Figure~\ref{fig:transfo_vltg} (left) illustrates 
the induced time-periodic voltages in the two coils, calculated with the time step $\delta T$ 
over the period $[0,0.02]$ s.

Linearization of the periodic system \eqref{eq:PP_PC_sys_explicit} 
is performed using the fixed point 
iteration \eqref{eq:add_splitting} with constant matrix
\begin{equation}
\setlength{\arraycolsep}{2.5pt}
\label{Hmatr_Klin}
\mathbf{H}:=
\begin{bmatrix}
\mathbf{C}+\bar{\mathbf{K}} &    &  & {-\mathbf{C}} \\
-\mathbf{C} & \mathbf{C}+\bar{\mathbf{K}} &  &  \\
&  \ddots   &  \ddots & \\
&     &  -\mathbf{C}    & \mathbf{C}+\bar{\mathbf{K}}
\end{bmatrix}
\end{equation}
where $\mathbf{C}={1}/{\Delta T}\cdot\mathbf{M}$ and 
$\bar{\mathbf{K}}$ is obtained for the fixed reluctivity in the core  
$\bar{\nu}=\SI{\nu_0\cdot 10^{-3}}{\metre\per\henry}$.  
An analogous linearization 
is applied also to the TP system \eqref{eq:sys_explicit_nlin}, 
using the matrix $\mathbf{H}$ from \eqref{Hmatr_Klin} with 
$\mathbf{C}={1}/{\delta T}\cdot\mathbf{M}.$ 
We now describe the 
performance of PP-PC MH and TP MH in terms of 
the effective and total numbers of the solved linear systems, as well as in terms of the wall clock time and compare it to the sequential time stepping. 



The standard sequential time integration with step size $\delta T=10^{-5}$ s required solution over $k^*=21$ periods until the steady state was reached 
in terms of the error \eqref{eq:error_tol_seq}. This corresponds to solution 
of $84\,000$ linear algebraic systems. 
It its turn, within the linearization of the form \eqref{eq:add_splitting} applied to nonlinear 
TP MH $600$ effective linear solutions were calculated during $6$ iterations. It therefore delivered the steady-state solution \IKRrr{solving $140$ times 
less effective linear systems} than the classical sequential approach. \IKRrr{In turn,} PP-PC MH converged 
\IKRrr{in two iterations, each requiring $6$ inner fixed point iterations \eqref{eq:add_splitting}.}  
It effectively needed $424$ linear systems solutions, \IKRrr{which is $198$ less} than with the sequential time stepping. \IKRrr{The total number of solved linear systems is equal to $6\times 2\,000=12\,000$ for TP MH and $8\,468$ for PP-PC MH. The corresponding percentages of the number of total and effective system solutions calculated by PP-PC MH and TP MH with regard to the sequential time stepping are illustrated Figure~\ref{fig:transfo_vltg} on the right.}

\IKRrr{Finally, measurement of the wall clock time showed that the TP MH calculation lasted slightly more than $9$ hours, while PP-PC MH only needed $2.6$ hours. 
The sequential computation over one period took almost $22$ hours, which over $21$ periods would extend to about $19$ days. The estimated speed up factors of TP MH and PP-PC MH are therefore about $51$ and $176$, respectively. We point out the considerable difference in the factors when comparing the effective number of solved linear systems and when comparing the actual execution time of the TP MH ($140$ vs. $51$) to the sequential calculation. This is due to the fact that the solution of the linear systems took about $60\%$ of the total computational time (while the DFT required $30\%$). In contrast to this, the measure based on the number of solved effective linear systems worked much better for estimating the computational effort of PP-PC MH, since the linear system solutions took the major part of the execution time.} 

\section{Conclusions}	\label{section:conclusions}
This paper presents an iterative parallel-in-time solution approach 
for nonlinear time-periodic problems combined with the MH 
coarse grid correction. The transformation into frequency domain 
introduces an additional parallelizability on the coarse grid, since the 
frequency components become decoupled. The proposed {simplified} Newton {method} 
with a special choice of the initial approximation keeps the nonlinearity 
constant over the Newton iterations and among the time instants,  
thereby letting the MH frequency domain approach be applicable 
at each iteration. Convergence of the iterative algorithm has been 
analyzed \IKR{in detail for a 1D model problem} and derived from the assumptions imposed on the nonlinearity. 
Performance of the method is illustrated for the time-periodic eddy current problem 
for both linear and nonlinear 2D models of a coaxial cable.   
Superiority of the MH solver over 
several existing (parallel-in-time) time-domain approaches is shown based on 
the computational costs calculated in terms of the number of required effective \IKRrr{and total} linear system solutions, \IKRrr{as well as when measuring the wall clock time}. 
\IKR{Numerical results demonstrate that the frequency domain transformation applied directly to the 
fine time-periodic system (within TP MH) deliver optimal scaling for the linear problem, 
while the MH correction performed at the simplified Newton iteration for the coarse PP-PC system performed better 
than TP MH in the nonlinear case for the majority of the considered settings (except $N=5$ CPUs).} 
\IKRr{Finally, the considered solution approaches with MH correction were applied to the steady-state analysis of a 3D
transformer model. It was observed that compared to the classical sequential time integration TP MH had 
\IKRrr{a speedup up to factor $51$} when using $N=20$ cores. A better result was obtained with the PP-PC MH algorithm, 
which was able to calculate the steady-state solution \IKRrr{using $176$} times less computational \IKRrr{time} than the sequential approach.}

\section*{Acknowledgements}
The authors would like to thank Herbert De Gersem for the multiple fruitful discussions 
on the multi-harmonic solution approach. 

Funding: This research was supported by the Excellence 
Initiative of the German Federal and State Governments and the Graduate School of 
Computational Engineering at Technische Universit\"at Darmstadt, as well as by DFG 
grant SCHO1562/1-2 and BMBF grant 05M2018RDA (PASIROM).


	\bibliographystyle{siamplain}
	

\end{document}